\documentclass[a4paper,reqno]{amsart}
\numberwithin{equation}{section}
\def\Cbbd{\mathbb{C}}

\def\Nbbd{\mathbb{N}}
\def\Qbbd{\mathbb{Q}}
\def\Rbbd{\mathbb{R}}

\def\Zbbd{\mathbb{Z}}
\def\Acal{\mathcal{A}}

\def\Fcal{\mathcal{F}}
\def\Hcal{\mathcal{H}}

\def\Pcal{\mathcal{P}}
\def\Qcal{\mathcal{Q}}

\def\Scal{\mathcal{S}}
\def\a{\alpha}
\def\d{\delta}
\def\l{\lambda}

\def\res{\mathop{\hbox{\rm res}}}
\def\id{\mathop{\hbox{\rm id}}\nolimits}
\def\abs#1{\left|#1\right|}
\def\ds{\displaystyle}

\def\endproof{\hfill\rule{2mm}{2mm}}
\def\?{(?)\marginpar{|?}}

\newtheorem{theo}{Theorem}[section]
\newtheorem{prop}{Proposition}[section]
\newtheorem{lemma}{Lemma}[section]
\newtheorem{remark}{Remark}[section]
\begin{document}
\title[Jack polynomials]%
  {$Q$-operator and factorised separation chain for Jack polynomials}
\author{Vadim B. Kuznetsov}
  \address{Department of Applied Mathematics,
          University of Leeds,
          Leeds LS2 9JT, UK}
  \email{V.B.Kuznetsov@leeds.ac.uk}
\author{Vladimir V. Mangazeev}
  \address{Centre for Mathematics and Its Applications, Mathematical Science
  Institute, Australian National University, Canberra, ACT 0200, Australia}
  \email{vladimir@maths.anu.edu.au}
\author{Evgeny K. Sklyanin}
\address{Department of Mathematics, University of York,
York YO10 5DD, UK}
\email{eks2@york.ac.uk }
\begin{abstract}
Applying Baxter's method of the $Q$-operator to the
set of Sekigu\-chi's commuting partial differential operators
we show that Jack polynomials
$P_\l^{(1/g)}(x_1,\dots,x_n)$ are eigenfunctions of
a one-parameter family of integral operators $Q_z$.
The operators $Q_z$ are expressed in terms of
the Dirichlet-Liouville $n$-dimensional beta integral.
{}From a composition of $n$ operators $Q_{z_k}$ we construct
an integral operator $\Scal_n$ factorising Jack polynomials into
products of hypergeometric polynomials of one variable.
The operator $\Scal_n$ admits a factorisation described in terms
of restricted Jack polynomials
$P_\l^{(1/g)}(x_1,\dots,x_k,1,\dots,\allowbreak 1)$.
Using the operator $Q_z$ for $z=0$ we give a simple derivation of
a previously known integral representation for Jack polynomials.
\end{abstract}
\noindent
{\it Dedicated to Tom H. Koornwinder on the occasion of his 60th birthday}
\begin{flushright}
to appear in {\it Indagationes Mathematicae}\\
math.CA/0306242
\end{flushright}
\maketitle
\tableofcontents

\section{Introduction}\label{intro}
\noindent
The method of $Q$-operator has been introduced by Rodney
Baxter \cite{Bax72,Bax82}
who used it to solve the eight-vertex model. Baxter's original $Q$-operator
was a finite matrix, which commuted with the transfer-matrix
of the model. Later on, Pasquier and Gaudin \cite{PG92} discovered
a connection between the classical B\"acklund transformation and a $Q$-operator
for the periodic Toda lattice. Their $Q$-operator was a certain
integral operator with a simple kernel. Recent developments
transformed the method into a unified approach used for solving
quantum integrable systems, with several families having been
extensively treated (see, e.g., \cite{BLZ,KS5,Der,KSS98,Pro,Skl53}).

Let us take a quantum integrable system defined by $n$ {\it commuting}
linear partial differential (or difference) operators $H_i$ in $n$ variables,
whose common eigenfunctions $\Psi_\l({\bf x})\equiv\Psi_\l(x_1,\dots,x_n)$,
\begin{equation}
H_i\Psi_\l({\bf x})=h_i(\l)\Psi_\l({\bf x}), \label{qop1}
\end{equation}
form a basis in a Hilbert space $\Hcal$.  The multi-index
$\l=(\l_1,\dots,\l_n)$ is a set of quantum numbers labelling
spectrum $h_i(\l)$ and eigenfunctions $\Psi_\l(\mbox{\bf x})$.

By definition, a $Q$-operator $Q_z$ depends on a parameter
$z\in\Cbbd$ and satisfies two commutativity properties
\begin{eqnarray}
{[Q_{z_1},Q_{z_2}]}&=&0\qquad\forall z_1,z_2\in\Cbbd,\label{qop2}\\
{[Q_{z},H_i]}&=&0\qquad \forall z\in \Cbbd,\quad \forall i=1,\dots,n,
\label{qop3}
\end{eqnarray}
which imply that $Q_z$ can be diagonalised by the same basis functions
$\Psi_\l(\mbox{\bf x})$:
\begin{equation}
[Q_z\Psi_\l]({\bf y})=q_\l(z)\Psi_\l({\bf y}).
\label{qop4}
\end{equation}

The third and the most important (as well as the characterising)
property of a $Q$-operator is that its eigenvalues $q_\l(z)$ satisfy
a linear {\it ordinary} differential (or difference) equation with
respect to the parameter $z$,
\begin{equation}
W\left(z,\frac{d}{dz};\{h_i(\l)\}\right)q_\l(z)=0,
\label{qop5}
\end{equation}
whose coefficients, apart from being functions in $z$, depend
on the eigenvalues $h_i(\l)$ of the commuting operators $H_i$.
Equation \eqref{qop5} is called {\it Baxter's equation\,} or, alternatively,
the {\it separation equation}.

An operator $Q_z$ for a quantum integrable system $\{H_i\}$
is a quantization of a B\"acklund transformation
for the corresponding classical integrable system.
In the quasi-classical limit, its kernel turns into the generating function
of the B\"acklund transformation, viewed as a canonical transform
\cite{PG92,KS5}.

A {\it separating operator} $\mathcal{S}_n$, by definition, factorises
the basis functions $\Psi_\l(\mbox{\bf x})$ into functions of {\it one}
variable:
\begin{equation}
\Scal_n: \Psi_\l(x_1,\dots,x_n) \mapsto c_\l\prod_{i=1}^n
\phi^{(i)}_\l(z_i)
\label{fact6}
\end{equation}
$c_\lambda$ being some normalization constants. The partial
functions $\phi^{(i)}_\l(z_i)$ are called {\it separated functions}
and usually satisfy a differential or difference equation like \eqref{qop5}.
Notice that, in the quasi-classical limit, the operator $\mathcal{S}_n$
turns into a separating canonical transformation for the corresponding
Liouville integrable system.

In \cite{KS5}, two of the authors noticed that any $Q$-operator
gives rise to a family of separating operators for an integrable system.
Indeed, given a $Q$-operator $Q_z$ and an arbitrary linear functional $\rho$
on $\Hcal$, one constructs a product of $n$ such operators
\begin{equation}
Q_{\bf z}=Q_{z_1}\cdots Q_{z_n}
\label{extrafact}
\end{equation}
with the kernel $Q_{{\bf z}}({\bf y}|{\bf x})$ as the convolution
\begin{equation}
Q_{{\bf z}}({\bf y}|{\bf x})=\int d{\bf t_1}\dots\int d{\bf t_{n-1}}\>
Q_{z_1}({\bf y}|{\bf t_1})Q_{z_2}({\bf t_1}|{\bf t_2})\dots
Q_{z_n}({\bf t_{n-1}}|{\bf x}),
\label{fact7}
\end{equation}
and defines an operator $\mathcal{S}_n^{(\rho)}=\rho Q_{\bf z}$ with the kernel
\begin{equation}
\Scal_n^{(\rho)}({\bf z}|{\bf x})
=\int d{\bf y} \>\rho({\bf y}) Q_{{\bf z}}({\bf y}|{\bf x}),
\label{efact8}
\end{equation}
where $\rho({\bf y})$ is the generalised function corresponding to $\rho$.
It then follows from \eqref{qop4} and
\eqref{fact6}--\eqref{efact8} that $\mathcal{S}_n^{(\rho)}$ is a family
of separating operators (parameterised by the functional $\rho$) with
\begin{equation}
\phi_\l^{(i)}(z_i)=q_\l(z_i)\qquad \mbox{\rm and}\qquad
c_\l=\int d{\bf y} \>\rho({\bf y})\Psi_\l({\bf y}).
\label{fact8}
\end{equation}

A clever choice of the `projection function' $\rho({\bf y})$
can, in principle, simplify the structure of the integral operator
$\mathcal{S}_n^{(\rho)}$.
Our main motivation is to study suitable classes of $Q$-operators
and functions $\rho({\bf y})$ for a variety of quantum integrable models
in order to construct
simplest separating operators $\mathcal{S}_n$ \eqref{fact6},
which factorise special functions $\Psi_\l(x_1,\dots,x_n)$.

In the present paper, we apply the above ideas to the Calogero-Sutherland
integrable model. In our case of study
$\mathcal{H}$ is 
an appropriate Hilbert space closure of
the space $\Cbbd[x_1,\dots,x_n]^{S_n}$
of symmetric polynomials in $n$ variables, $\Psi_\l(x_1,\dots,x_n)$
are Jack polynomials and $\l=(\l_1,\dots,\l_n)$ is a partition.

The problem of separation of variables for the Calogero-Sutherland
model and for its $q$-analogue has been addressed by the authors
in a series of publications since 1994.
Many of the results presented below generalise those obtained from
the studies of the few particle cases, $n=2,3,4$, in
\cite{Skl38,KS4,KS55,KNS,KS96,KS99,Man,Man2}.

The paper is organised as follows. Section \ref{jack} fixes notation and
collects the necessary background information on Jack polynomials.
In Section \ref{dirichlet} we study the properties of 
Dirichlet-Liouville multidimensional beta-integrals, which constitute our
main analytical tool. All integral operators that we construct are expressed
in terms of these integrals. In Section \ref{q} we introduce a Q-operator
$Q_z$ for the Calogero-Sutherland model and derive its properties.
In particular, we prove that $Q_z$ is diagonal in the basis of Jack 
polynomials. Thus, the eigenvalue problem $Q_z\Psi_\l=q_\lambda(z)\Psi_\lambda$
can be considered as an integral equation for Jack polynomials.
In Section \ref{eigenvalues} we obtain the eigenvalues 
$q_\lambda(z)$ of $Q_z$ and give for them several explicit expressions,
including the one in terms of the generalised hypergeometric series.
In section \ref{restricted} we construct the separating operator
$\mathcal{S}_n^{(\rho)}$ and show that for a special choice of 
the functional $\rho=\rho_0$ it splits into a product
of $n$ operators 
 $\Acal_i$,
$i=1,\dots,n$,
\begin{equation}
\Scal_n^{(\rho_0)}=\Acal_1\cdots \Acal_n,
\label{fact10a}
\end{equation}
such that each operator $\Acal_k$ acts only on $k$ variables.
Such a factorisation is an example of the {\it factorised separation chain}
introduced in \cite{Kuzn} for integrable models of $sl(2)$ class.

In Section \ref{intrepr} we describe another application of the $Q$-operator
$Q_z$.
Setting $z=0$ we show that the operator $Q_0$ leads to
an integral relation between Jack polynomials in
$n$ and $n-1$ variables $\{x_i\}$. Iterating
this relation we produce an integral representation
for Jack polynomials, which is equivalent
to the integral representation due to Okounkov and Olshanski \cite{Oko} 
who found it as a limit from a formula involving
shifted Jack polynomials.

\section{Jack polynomials}\label{jack}
\noindent
This section serves an introductory purpose as it fixes notation
and collects the necessary background information on Jack polynomials.

Jack polynomials $P_{\l}^{(\alpha)}(x_1,\dots,x_n)$
\cite{Jack,Jack2,Jack3,Koorn,Stan,Macd}
are orthogonal symmetric polynomials, discovered by Henry Jack
in 1970 as a one-parameter extension of the Schur functions $s_\l(x_1,\dots,x_n)$.
In the case $\alpha=2$, they reduce to zonal polynomials constructed
by James in 1960 (\cite{James}).
We shall use the parameter $g=1/\alpha$, instead of $\alpha$ as in \cite{Macd}, i.e.
use the notation $P_\l^{(1/g)}(\mbox{\bf x})\equiv P_{\l}^{(1/g)}(x_1,\dots,x_n)$.

Let $\l=(\l_1,\l_2,\dots,\l_n)\in\Nbbd^n$, $\l_1\ge\l_2\ge \ldots\ge\l_n\ge 0$,
be a partition of an arbitrary weight $|\l|$,
\begin{equation}
|\l|=\sum_{i=1}^n\l_i.
\label{jack0}
\end{equation}
Hereafter, we shall use the weight notation $|\,.\,|$ as in (\ref{jack0}) for any
finite sum of indexed  variables.
The dominance partial ordering $\preceq$ for two partitions $\mu$ and $\l$
is defined as follows:
\begin{equation}
 \mu\preceq \l \quad \Longleftrightarrow \quad
\Bigl\{ |{\mu}|=|{\l}|\,; \quad
\sum_{j=1}^k \mu_j\leq\sum_{j=1}^k \l_j\,, \quad
k=1,\dots,n-1\Bigr\}\,.
\label{jack5}
\end{equation}

Usually (cf \cite{Macd}), the {\it length}\, $l(\l)$ of a partition $\l$
is the number of non-zero {\it parts} $\l_i$.   We shall consider
the partitions $\l=(\l_1,\l_2,\dots,\l_n)$ that are finite sequences
of exact length $n$, {\it including}\, any zero parts. For instance, $(6,6,2,0)$
and $(6,0,0,0)$ are both partitions of the length $4$, while $(6,6,2)$ and
$(6)$ have lengths $3$ and $1$, respectively. The reason for this deviation
{} from the standard definition is because we shall study Jack
polynomials with the {\it finite} \,number of variables and characterise them as
eigenfunctions of the complete set of $n$ commuting
partial differential operators $\{H_i\}$ in $n$ variables.
In this picture, the parts $\l_i$ are the labels that parameterise
the spectrum and the eigenfunctions of the operators
$\{H_i\}$. The number of parts, therefore,
will always be equal to the number of the operators $H_i$ (or to the
number of the variables $x_i$). 

In the main text, for any partition $\l=(\l_1,\dots,\l_n)$ of the length $n$,
we shall often use the notation $\l_{i,j}=\l_i-\l_j$,  in particular
$\l_{i,i+1}=\l_i-\l_{i+1}$, $i=1,\dots,n$, $\l_{n+1}\equiv 0$.

In the theory of symmetric polynomials, the number $n$ of the variables $x_i$ is
often
irrelevant because many results have been obtained for the case when this number is
assumed to be large enough, for instance $n\geq |\l|$, for the
{\it stabilised} polynomials \cite{Stan}. In this paper, the number $n$
of the variables $x_i$ is a fixed finite number, the partition $\l$ is of the
length $n$ and has an
{\it arbitrary} weight $|\l|$, which means that we shall mostly deal with
the {\it non-stabilised} Jack polynomials.

Let $\Cbbd(a,b,c,\ldots)$ be the field of rational functions in indeterminates
$a$, $b$, $c,\ldots$ over the field $\Cbbd$ of complex numbers,
$\Cbbd(a,b,c,\ldots)[x,y,z,\ldots]$ be the ring of polynomials
in the variables $x$, $y$, $z,\ldots$ with coefficients from
$\Cbbd(a,b,c,\ldots)$, and $\Cbbd[x_1,\dots,\allowbreak x_n]^{S_n}$ be the
subring of symmetric polynomials.

We need two standard bases in $\Cbbd[\mathbf{x}]^{S_n}$ labelled
by partitions $\l$.
For each $\mbox{\bf a}=(a_1,\dots,a_n)\in \Nbbd^n$ let us denote by  ${\bf x}^{\bf
a}$ the monomial
\begin{equation}
{\bf x}^{\bf a}=x_1^{a_1}\cdots x_n^{a_n}.
\end{equation}
The {\it monomial symmetric functions} $m_\l(\mathbf{x})$ are defined by
\begin{equation}
 m_\l(\mathbf{x})=\sum {\bf x}^\nu,
\label{def-m}
\end{equation}
where the sum is taken over all distinct permutations $\nu$ of $\l$.
The $m_\l(\mathbf{x})$ form a basis in $\Cbbd[\mathbf{x}]^{S_n}$.

For each $r=0,\dots,n$ the $r$th {\it elementary symmetric function}
$e_r(\mathbf{x})$ is the sum of all products of $r$ distinct variables
$x_i$, so that $e_0({\bf x})=1$ and
\begin{equation}
e_r(\mathbf{x})=\sum_{1\leq i_1<\ldots<i_r\leq n}x_{i_1}x_{i_2}\cdots
x_{i_r}=m_{(1^r0^{n-r})}({\bf x})
\label{qop12}
\end{equation}
for $r=1,\dots,n$.
For each partition $\l=(\l_1,\dots,\l_n)$ define
the polynomials $E_{\l}(\mathbf{x})$ as
\begin{equation}
E_{\l}(\mathbf{x})
  =e_1^{\l_1-\l_2}(\mathbf{x})e_2^{\l_2-\l_3}(\mathbf{x})
\cdots e_n^{\l_n}(\mathbf{x}).
\label{qop37}
\end{equation}
The polynomials $E_{\l}(\mathbf{x})$ also form a basis in $\Cbbd[\mathbf{x}]^{S_n}$
and  the transition matrix between $E_{\l}$ and $m_{\l}$ is triangular
with respect to the dominance partial ordering (see \cite{Macd}, Ch.\ 1, Sec.\ 2):
\begin{equation}
E_{\l}(\mathbf{x})=\sum_{\mu\preceq\l}a_{\l\mu}m_\mu(\mathbf{x}),
\label{qop40}
\end{equation}
where  $a_{\l\mu}\in\Nbbd$ and $a_{\l\l}=1$.

Jack polynomials $P_\l=P_\l^{(1/g)}(\mathbf{x})
\in\Cbbd(g)[\mathbf{x}]^{S_n}$
are homogeneous of degree $\abs{\lambda}$:
\begin{equation}
P_\l^{(1/g)}(x_1,\dots,x_n)=x_n^{|\l|}P_\l^{(1/g)}(x_1/x_n,
\dots,x_{n-1}/x_n,1).
\label{jack14}
\end{equation}
They form a basis in $\Cbbd(g)[\mathbf{x}]^{S_n}$. The transition
matrix between the bases $P_\l$ and $m_\l$
is triangular with respect to the dominance partial ordering:
\begin{equation}
P_\l=\sum_{\mu\preceq\l}v_{\l\mu}(g)m_{\mu},
\qquad v_{\l\mu}(g)\in\Cbbd(g),
\quad v_{\l\l}=1.
\label{jack4}
\end{equation}
The three bases, $m_\lambda$, $E_\lambda$ and $P_\lambda$, are thus
mutually triangular.

The polynomials $P_\l^{(1/g)}(\mathbf{x})$ diagonalise the
Sekiguchi partial differential operators \cite{Sek} given by the generating
function $D(u;g)$:
\begin{equation}
D(u;g)=\frac{(-1)^{n(n-1)/2}}
{\Delta_n({\bf x})}\,\det\biggl[x_i^{n-j}\Bigl(x_i\frac{\partial}
{\partial x_i}+(n-j)g+u\Bigr)\biggr]_{i,j=1}^n,\label{jack6}
\end{equation}
\begin{equation}
D(u;g)\;P_\l^{(1/g)}(\mathbf{x})=\left[\prod_{i=1}^n(\l_i+(n-i)g+u)\right]\,P_\l^{(1/
g)}(\mathbf{x}),
\label{jack7}
\end{equation}
where
\begin{equation}
\Delta_n({\bf x})=\det\left[x_j^{i-1}\right]_{i,j=1}^n=\prod_{1\le i<j\le
n}(x_j-x_i)\label{jack8}
\end{equation}
is the Vandermonde determinant. Expansion \eqref{jack4} and the fact that $P_\l$
are the eigenfunctions of the Sekiguchi operators characterise
Jack polynomials uniquely.

Generating function $D(u;g)$ is an operator polynomial in $u$
whose coefficients give $n$ commuting partial differential
operators $H_i$, $i=1,\dots,n$:
\begin{equation}
D(u;g)=u^n+\sum_{i=1}^n u^{n-i}H_i,\qquad
H_1=\frac{n(n-1)}{2}\,g+\sum_{i=1}^nx_i\,\frac{\partial}{\partial x_i}\,.
\label{jack10}
\end{equation}
Define the second order partial differential operator $H$ as
\begin{equation}
H=\frac{n(n-1)(n-2)}{6}\,g^2+H_1[H_1-(n-1)g]-2H_2\label{jack11}
\end{equation}
then
\begin{equation}
H=\sum_{i=1}^n\Bigl(x_i\,\frac{\partial}{\partial x_i}\,\Bigr)^2+
g\sum_{i<j}\frac{x_i+x_j}{x_i-x_j}\,\Bigl(x_i\,\frac{\partial}{\partial x_i}-
x_j\,\frac{\partial}{\partial x_j}\,\Bigr).\label{jack12}
\end{equation}
It follows from (\ref{jack7}) that
\begin{equation}
H\,P_\l^{(1/g)}(\mathbf{x})=h(\l)\, P_\l^{(1/g)}(\mathbf{x}),\qquad
h(\l)=\sum_{i=1}^n
\l_i[\l_i+g(n+1-2i)].
\label{jack13}
\end{equation}

Jack polynomials were applied in \cite{Cal, Suth, OP}
to describe the excited states in the Calogero-Sutherland model.
Let us briefly remind some facts about this model.
It describes a system of $n$ quantum particles on a circle
 with coordinates $0\le q_i\le \pi$, $i=1,\dots,n$.
The quantum Hamiltonian and the operator of the total momentum are
 \begin{equation}
 H_{CS}=-\frac{1}{2}\sum_{i=1}^n\frac{\partial^2}{\partial q_i^2}+\sum_{i<j}
 \frac{g(g-1)}{\sin^2(q_i-q_j)}\,,\qquad P_{CS}=-i\sum_{j=1}^n\frac{\partial}{\partial q_j}\,.
\label{jack15}
 \end{equation}
The space of quantum states of the model is the complex Hilbert space
$L^2(T^n/S_n)$
of functions $\Psi(\mathbf{q})$ on the torus $T^n={\Rbbd^n}/\pi{\Zbbd^n}\ni
\mathbf{q}$, symmetric with respect to permutations of the coordinates
$q_i$. The scalar product in $L^2(T^n/S_n)$
is defined as follows:
\begin{equation}
\langle\Psi,\Phi\rangle=\frac{1}{n!}
  \int_{0}^\pi dq_1\cdots\int_0^\pi dq_n\,
\overline{\Psi(\mathbf{q})}\,\Phi(\mathbf{q}).
\label{jack16}
\end{equation}

The differential operators \eqref{jack15} are formally Hermitian
with respect to the scalar product \eqref{jack16} when the parameter $g$
is real.
Most results of the present paper are established, however, for $g>0$,
which we shall assume hereafter.

The ground state of the model is given by the function
\begin{equation}
\omega(\mathbf{q})=\biggl|\,\prod_{i<j}\sin(q_i-q_j)\biggr|^g
\label{jack17}
\end{equation}
with the ground state energy $E_0=\frac{1}{6}\,g^2(n^3-n)$.
The gauge transformation of the
Hamiltonian $H_{CS}$ (\ref{jack15}) with the function $\omega(\mathbf{q})$
gives
\begin{equation}
\omega^{-1}(\mathbf{q})\circ H_{CS}\circ\omega(\mathbf{q})=\frac{1}{2}\,H+E_0,
\label{jack19}
\end{equation}
where the operator $H$ is defined by (\ref{jack12}) and the change
of variables from $\mathbf{q}$ to $\mathbf{x}=\varepsilon(\mathbf{q})$ is given by the map
\begin{equation}
   \varepsilon:(q_1,\dots,q_n)\mapsto(e^{2iq_1},\dots,e^{2iq_n}).
\label{jack18a}
\end{equation}
The eigenfunctions of the $H_{CS}$ can therefore
be chosen as the product of the ground state function and
a Jack polynomial labelled by a partition $\l$:
\begin{eqnarray}
\Psi_\l(\mathbf{q})=\omega(\mathbf{q})P_\l^{(1/g)}(\varepsilon(\mathbf{q})).
\label{jack18}
\end{eqnarray}
The Hermitian form (\ref{jack16}) corresponds to the following
scalar product for any two symmetric polynomials
$f(\mathbf{x})$ and $p\,(\mathbf{x})$ from $\Cbbd[\mathbf{x}]^{S_n}$:
\begin{equation}
 \langle f,p\,\rangle=\int_{T_\mathbf{v}} d\mathbf{v}\;
\left[\omega(\mathbf{v})\right]^2
   \, \overline{f(\varepsilon(\mathbf{v}))}\;p\,(\varepsilon(\mathbf{v})),
\label{jack20}
\end{equation}
where the integration domain, by symmetry, is
\begin{equation}
 T_\mathbf{v}=\{\mathbf{v}\in\Rbbd^n/\pi\Zbbd^n \mid
     v_j\in(v_{j-1},v_{j+1}), \quad j+n\equiv j\}.
\label{defTv}
\end{equation}
The intervals $(v_{j-1},v_{j+1})$ above are
understood as arcs on the circle $\Rbbd/\pi\Zbbd$
such that $v_{j-1}$ precedes $v_{j+1}$
with respect to the orientation on $\Rbbd$.

For $g>0$, Jack polynomials
$P_\l^{(1/g)}(\mathbf{x})$ are orthogonal with respect to the
scalar product \eqref{jack20}.
Note that
${\overline{P_\l^{(1/g)}(\mathbf{x})}}
=P_\l^{(1/g)}(\overline{x_1},\allowbreak\dots,\overline{x
_n}\,)$.

\section{Dirichlet-Liouville integral and its modifications}\label{dirichlet}
\noindent
The Dirichlet integrals and their extension due to Liouville (see, e.g., \cite{AAR})
are the most straight\-forward multidimensional generalisations of the beta-integral.
They are used in computing volumes and surface areas
of multidimensional ellipsoids. A special form of these integrals
serves as the main building tool for our construction
of the $Q$-operator for Jack polynomials. In what follows,
the powers $\xi^\a$ are functions of a positive argument $\xi>0$
and are understood as $\xi^\a=\exp(\a\ln\xi)$, with
$\ln\xi\in\Rbbd$.

\begin{theo}[Dirichlet-Liouville \cite{AAR}]\label{liov1}
Define the domain $V_{\boldsymbol{\xi}}\subset \Rbbd^n$ as
\begin{equation}
V_{\boldsymbol{\xi}}=\{\boldsymbol{\xi}\in\Rbbd^n \mid \xi_i>0,\, i=1,\dots,n\}.
\label{def-Vxi}
\end{equation}
Then for $\Re\a_i>0$,
\begin{equation}
\int_{V_{\boldsymbol{\xi}}} d\boldsymbol{\xi}\;\,\d(\abs{\boldsymbol{\xi}}-1)\;\,
\xi_1^{\a_1-1}\cdots\,\xi_n^{\a_n-1}\,
=\frac{\Gamma(\a_1)\cdots\Gamma(\a_n)}{\Gamma(\abs{\boldsymbol{\a}})}\,,
\label{liov2}
\end{equation}
where $\boldsymbol{\a}=(\a_1,\dots,\a_n)$, 
$|\boldsymbol{\a}|=\a_1+\cdots+\a_n$,
$d\boldsymbol{\xi}=d\xi_1\cdots d\xi_n$,
and $\d$ is the Dirac delta-function.
\end{theo}
The proof can be found in \cite{AAR}.
We use two modifications of the integral \eqref{liov2} given in Theorems
\ref{dirichlet1} and \ref{dirichlet2} below.

\begin{theo}\label{dirichlet1}
Let $\mathbf{y}=(y_1,\dots,y_n)$ be a fixed
set of real positive parameters such that
\begin{equation}
0< y_1< y_2< \ldots< y_n<\infty.
\label{liov2a}
\end{equation}
Define the domain $\Omega_{\mathbf{x}}\subset \Rbbd^n$ by the inequalities
\begin{equation}
\Omega_{\mathbf{x}}=
\{\mathbf{x}\in \Rbbd^n \mid 0< y_1< x_1< y_2<\ldots< x_{n-1}< y_n< x_n<\infty\}.
\label{liov6}
\end{equation}
Let $\check{x}\equiv e_n(\mathbf{x})=x_1\cdots x_n$
and \,$\check{y}\equiv e_n(\mathbf{y})= y_1\cdots y_n$.
Then for $\Re\a_i>0$ and $z>1$,
\begin{equation}
\int_{\Omega_{\mathbf{x}}}\, d\mathbf{x}\;
\delta(\check{x} - z\check{y})\,
\Delta_n(\mathbf{x})\,
 \prod_{i=1}^n\xi_i^{\a_i-1}
=(z-1)^{\abs{{\boldsymbol{\a}}}-1}\Delta_n(\mathbf{y})\,
\frac{\prod_{i=1}^ny_i^{\a_i-1}\Gamma(\a_i)}{\Gamma(\abs{{\boldsymbol{\a}}})}\,,
\label{liov13}
\end{equation}
where $\Delta_n$ is the Vandermonde determinant \eqref{jack8} and
$\xi_i$ are defined as follows:
\begin{equation}
\xi_i=\frac{\prod_{k=1}^n(x_k-y_i)}{\prod_{k\neq i}(y_k-y_i)}\,.
\label{liov5}
\end{equation}
\end{theo}

The other modification of the Dirichlet-Liouville integral
is a trigonometric variant of \eqref{liov13}.
\begin{theo}\label{dirichlet2}
Given a fixed set of real parameters $\mathbf{v}=(v_1,\dots,v_n)$ such that
\begin{equation}
 v_1< v_2<\ldots< v_n< \pi+v_1,
\label{liov16}
\end{equation}
define the domain $\Omega_{\mathbf{u}}^T\subset \Rbbd^n$ by the inequalities
\begin{equation}
\Omega^T_{\mathbf{u}}=\{\mathbf{u}\in \Rbbd^n \mid
 v_1< u_1< v_2<\ldots< u_{n-1}< v_n< u_n<\pi+v_1\}.
\label{liov18}
\end{equation}
Then for $\gamma\in(0,\pi)$ and $\Re\a_i>0$,
\begin{multline}
\int_{\Omega^T_{\mathbf{u}}}d{\mathbf{u}}\,
\delta(\abs{\mathbf{u}}-\abs{\mathbf{v}}-\gamma)
\left[\prod_{i<j}\sin(u_j-u_i)\right]\,
\left[\prod_{i=1}^n(\xi^T_i)^{\a_i-1}\right] \\
=(\sin\gamma)^{\abs{\boldsymbol{\a}}-1}\,\left[\prod_{i<j}\sin(v_j-v_i)\right]\,
\frac{\prod_{i=1}^n\Gamma(\a_i)}{\Gamma(\abs{\boldsymbol{\a}})}\,,
\label{liov23}
\end{multline}
where
\begin{equation}
\xi^T_i=\frac{\prod_{k=1}^n\sin(u_k-v_i)}{\prod_{k\neq i}
\sin(v_k-v_i)}\,.
\label{liov17}
\end{equation}
\end{theo}

Both Theorems  are proved by
performing the respective changes of variables, from $\{\xi_i\}$ to $\{x_i\}$
\eqref{liov5} and from $\{\xi_i^T\}$ to $\{u_i\}$ \eqref{liov17}, in the integral
\eqref{liov2}.

{\bf Proof of Theorem \ref{dirichlet1}.}
The first step is to rescale the variables $\xi_i\mapsto \frac{\xi_i}{y_i(z-1)}$ in
the integral
\eqref{liov2} and to rewrite it in the form
\begin{equation}
\int_{V_{\boldsymbol{\xi}}} d{\boldsymbol{\xi}}\;
\d\left(\sum_{i=1}^n\frac{\xi_i}{y_i}-(z-1)\right)
\prod_{i=1}^n\xi_i^{\a_i-1}
=(z-1)^{\abs{\boldsymbol{\a}}-1}\,\frac{\prod_{i=1}^n
y_i^{\a_i}\Gamma(\a_i)}{\Gamma(\abs{\boldsymbol{\a}})}\,.
\label{liov3}
\end{equation}

The next step is to prove
that the relations \eqref{liov5} define a bijection of
$V_{\boldsymbol{\xi}}$ onto $\Omega_{\mathbf{x}}$. From the inequalities
in \eqref{liov6} it follows immediately that
$\Omega_{\mathbf{x}}$ is mapped into $V_{\boldsymbol{\xi}}$.
To prove the invertibility of this mapping consider the rational function
\begin{equation}
R(t)=\prod_{i=1}^n\frac{t-x_i}{t-y_i}\,.
\label{liov8}
\end{equation}
It is easy to see that $\res_{t=y_i} R(t)=-\xi_i$
and $\lim_{t\rightarrow\infty}R(t)=1$, and therefore
\begin{equation}
R(t)=1-\sum_{i=1}^n\frac{\xi_i}{t-y_i}.
\label{Rres}
\end{equation}
Given a set of positive $\xi_i>0$, the function $R(t)$ is restored by the
formula \eqref{Rres}, and the variables $x_i$ are then determined
uniquely as its zeroes, taken in the increasing
order. The standard continuity argument,
using $\xi_i>0$, shows that
all $n$ zeroes $x_j$ of $R(t)$ are real and interspersed with
the poles $y_j$ of $R(t)$  and $+\infty$ or, equivalently,
belong to the domain $\Omega_{\mathbf{x}}$.

The bijectivity of the mapping $\boldsymbol{\xi}\mapsto\mathbf{x}$
having been thus established, its Jacobian is as
follows:
\begin{equation}
J=\det\biggl[\frac{\partial\xi_j}{\partial x_i}\biggr]_{i,j=1}^n=
\det\biggl[\frac{\xi_j}{x_i-y_j}\biggr]_{i,j=1}^n=
\frac{\Delta_n(\mathbf{x})}{\Delta_n(\mathbf{y})}>0,
\label{liov9}
\end{equation}
which was evaluated with the help of Cauchy's determinantal formula \cite{Cauchy}
\begin{equation}
\det\left[   \frac{1}{1-a_ib_j}\right]_{i,j=1}^n=\frac{\Delta_n(\mathbf{a})\Delta_n(\mathbf{b
})}
{\prod_{i,j=1}^n(1-a_ib_j)}\,.
\label{liov10}
\end{equation}
Setting $t=0$ in \eqref{liov8} and \eqref{Rres} we get
$1+\sum_{i=1}^n\xi_i/y_i=\prod_{k=1}^n x_k/y_k$, and
hence,
\begin{equation}
\delta\left( \sum\frac{\xi_i}{y_i}-(z-1)\right)
=\left(\prod y_k\right) \;\delta\left(\prod x_k - z \prod y_k\right).
\label{liov12}
\end{equation}
Using \eqref{liov9} and \eqref{liov12} we can now fulfil the
change of variables $\boldsymbol{\xi}\mapsto\mathbf{x}$ in \eqref{liov3} and
obtain \eqref{liov13}.
\endproof

The proof of Theorem \ref{dirichlet2} is similar to the above, but
the structure of the mapping $\boldsymbol{\xi}^T\mapsto\mathbf{u}$ \eqref{liov17}
is, however, more involved compared to the rational case \eqref{liov5}.

\begin{prop}
Define the domains $V_{\boldsymbol{\xi}}^T,\Phi_{\mathbf{u}}^\pm\subset\Rbbd^n$ as
\begin{equation}
  V_{\boldsymbol{\xi}}^T=\{\boldsymbol{\xi}\in\Rbbd^n \mid
     \abs{\boldsymbol{\xi}}<1,\; \xi_i>0;\, i=1,\dots,n\}
\label{defVxiT}
\end{equation}
and
\begin{equation}
  \Phi_{\mathbf{u}}^-=\{\mathbf{u}\in\Omega_{\mathbf{u}}^T \mid
\abs{\mathbf{u}}<\abs{\mathbf{v}}+\frac{\pi}{2}\,\},
\qquad
  \Phi_{\mathbf{u}}^+=\{\mathbf{u}\in\Omega_{\mathbf{u}}^T \mid
\abs{\mathbf{u}}>\abs{\mathbf{v}}+\frac{\pi}{2}\,\}.
\label{defPhipm}
\end{equation}
Then the change of variables 
$\mathbf{u}\mapsto\boldsymbol{\xi}^T$ maps the sets
$\Phi_{\mathbf{u}}^{\pm}$ bijectively onto $V_{\boldsymbol{\xi}}^T$.
\end{prop}

{\bf Proof.} Consider the function $R(t)$
\begin{equation}
R(t)=\prod_{i=1}^n\frac{\sin(t-u_i)}{\sin(t-v_i)}.
\label{defRtrig}
\end{equation}
Noting that $\res_{t=v_i}R(t)=-\xi_i^T$, $R(t)$ can be
rewritten in the form
\begin{equation}
 R(t)=c-\sum_{i=1}^n \xi_i^T\cot(t-v_i).
\label{Rtrigexpand}
\end{equation}
Comparing the limits of the two expressions for $R(t)$ at
$t\rightarrow\pm i\infty$ one gets
\begin{equation}
 c=\cos(\abs{\mathbf{u}}-\abs{\mathbf{v}}), \qquad
 \abs{{\boldsymbol{\xi}}^T}=\sin\left(\abs{\mathbf{u}}-\abs{\mathbf{v}}\right).
\label{defc}
\end{equation}

Suppose that $\mathbf{u}\in\Omega_{\mathbf{u}}^T$, that is either
$\mathbf{u}\in\Phi_{\mathbf{u}}^-$ or $\mathbf{u}\in\Phi_{\mathbf{u}}^+$.
Then, from the inequalities \eqref{liov18} and the expressions \eqref{liov17}
for $\xi_i^T$, it follows that $\xi_i^T>0$. It also
follows from \eqref{defc} that $\abs{{\boldsymbol{\xi}}^T}\in(0,1)$, for
$0<\abs{\mathbf{u}}-\abs{\mathbf{v}}<\pi$. Thus,
both domains $\Phi_{\mathbf{u}}^\pm$ are  mapped into
$V_{\boldsymbol{\xi}}^T$.
Conversely, given $\boldsymbol{\xi}^T\in V_{\boldsymbol{\xi}}^T$, one obtains
{} from the relation $c^2+\abs{\boldsymbol{\xi}^T}^2=1$ (cf \eqref{defc})
two possible values for the constant $c$ in \eqref{Rtrigexpand}.
The inequalities in \eqref{defVxiT} ensure that
$c^2\in(0,1)$. Having chosen a solution for $c$, one reconstructs
$R(t)$ using the formula \eqref{Rtrigexpand}. The continuity argument,
based on the positivity of $\xi_i^T$, shows that the zeroes $u_i$ of
$R(t)$ belong to the domain $\Omega_{\mathbf{u}}^T$. The positive
$c\in(0,1)$ correspond  to
$\mathbf{u}\in\Phi_{\mathbf{u}}^-$ and
the negative $c\in(-1,0)$ to $\mathbf{u}\in\Phi_{\mathbf{u}}^+$.
\endproof

{\bf Proof of Theorem \ref{dirichlet2}} starts with rescaling the
variables $\xi_i\mapsto\frac{\xi_i}{\sin \gamma}$ in the
integral \eqref{liov2} which transforms it into the formula
\begin{equation}
\int_{V_{\boldsymbol{\xi}}}
\,d{\boldsymbol{\xi}}\;\d\left(\abs{\boldsymbol{\xi}}-\sin\gamma\right)\;
\xi_1^{\a_1-1}\cdots\xi_n^{\a_n-1}\,
=(\sin\gamma)^{\abs{\boldsymbol{\a}}-1}\,
\frac{\Gamma(\a_1)\cdots\Gamma(\a_n)}{\Gamma(\abs{\boldsymbol{\a}})}\,.
\label{liov2b}
\end{equation}
By the condition of the Theorem, $\gamma\in(0,\pi)$, so
$\sin\gamma\in(0,1)$, and therefore the integration
in \eqref{liov2b} can be restricted to $V_{\boldsymbol{\xi}}^T$.
The Jacobian of the mapping $\boldsymbol{\xi}^T\mapsto\mathbf{u}$,
\begin{equation}
J^T=\det\left[\frac{\partial\xi_j^T}{\partial
u_i}\right]_{i,j=1}^n=\;\det[\xi_j^T\cot(u_i-v_j)]_{i,j=1}^n,
\label{liov22}
\end{equation}
can be calculated with the help of a trigonometric limit of the
generalized (elliptic) Cauchy determinant (cf \cite{Frobenius,KNS}),
\begin{equation}
\det\left[\frac{\sigma(x_i-y_j+\zeta)}{\sigma(x_i -
y_j)\sigma(\zeta)}\right]_{i,j=1}^n =
\frac{\sigma(\Sigma+\zeta)}{\sigma(\zeta)}\;
\frac{ \prod_{k<l} \sigma(x_k - x_l)
\,\sigma(y_l - y_k) }{\prod_{k,l} \sigma(x_k - y_l) }\,,
\label{eq:cauchy} 
\end{equation}
where $\sigma(x)$ is the Weierstrass sigma-function
and $\Sigma\equiv \sum_i (x_i-y_i)$. Take
the limit of the infinite imaginary
period, $\omega^\prime\rightarrow\infty$, choose the real period
$\omega=\frac{\pi}{2}$ and put $\zeta=\frac{\pi}{2}$,
then $\sigma(x)\rightarrow\exp\left(\frac16\,x^2\right)\sin x$ and
\eqref{eq:cauchy} turns into the following trigonometric
determinantal formula:
\begin{equation}
\det\left[\cot(x_i-y_j)\right]_{i,j=1}^n =
 \cos(\Sigma )\;  \frac{ \prod_{k<l} \sin(x_k - x_l)
\,\sin(y_l - y_k) }{  \prod_{k,l} \sin(x_k - y_l) }\,.
\label{eq:cauchy2} 
\end{equation}
Now we can evaluate the Jacobian $J^T$ in \eqref{liov22}:
\begin{equation}
J^T=\cos(\abs{\mathbf{u}}-\abs{\mathbf{v}})
\prod_{i<j}\frac{\sin(u_j-u_i)}{\sin(v_j-v_i)}\,.
\label{liov22bis}
\end{equation}
Notice that $J^T$ is positive  for $\mathbf{u}\in\Phi_{\mathbf{u}}^-$ and
negative  for $\mathbf{u}\in\Phi_{\mathbf{u}}^+$.

Using \eqref{defc}, one transforms the $\d$-function in \eqref{liov2b} as
\begin{equation}
\d(\abs{\boldsymbol{\xi}}-\sin\gamma)=\frac{\d(\abs{\mathbf{u}}-\abs{\mathbf{v}}-\gamma)}{\cos\gamma}.
\label{deltatrig}
\end{equation}

Finally, combining \eqref{liov22bis} and \eqref{deltatrig} one fulfils the
change of variables $\boldsymbol{\xi}^T\mapsto\mathbf{u}$ in
\eqref{liov2b} and obtains \eqref{liov23}. When $\gamma\in(0,\pi/2)$,
the integration in \eqref{liov23} can be restricted to $\Phi_{\mathbf{u}}^-$, and
when $\gamma\in(\pi/2,\pi)$, the integration is over  $\Phi_{\mathbf{u}}^+$.
Formula \eqref{liov23} is valid for any $\gamma\in(0,1)$ by continuity.
\endproof

\section{$Q$-operator $Q_z$}\label{q}
\noindent
In this section we construct a $Q$-operator $Q_z$ for the Sekiguchi
partial differential operators $\{H_i\}$ \eqref{jack10}.
The operator $Q_z$ is described
as an integral operator with a kernel 
expressed in terms of elementary functions
and satisfies the defining conditions
\eqref{qop2}--\eqref{qop4}, proved in this section, and \eqref{qop5},
proved in the next section.
Thus, the eigenvalue problem $Q_zP_\l^{(1/g)}=q_\lambda(z)P_\l^{(1/g)}$
can be considered as an integral equation for Jack polynomials
$P_\l^{(1/g)}$.

For $g>0$, $z>1$ and for an ordered set of $n$ real positive parameters
$\mathbf{y}=(y_1,\dots,y_n)$,
\begin{equation}
0< y_1< y_2< \ldots< y_n<\infty,
\label{liov2abis}
\end{equation}
define the operator $Q_z$, acting on symmetric
polynomials $p\in\Cbbd[\mathbf{x}]^{S_n}$,
by the integral
\begin{equation}
  [Q_zp](\mathbf{y})=
  \frac{\Gamma(ng)}{\Gamma(g)^n}\,
  \frac{(z-1)^{1-ng}}{\ds\Delta_n({\bf y}\,)\prod_{i=1}^n y_i^{g-1}}\,
  \int_{\Omega_{\mathbf{x}}}\, d{\mathbf{x}}\;
\delta(\check{x} - z\check{y})\,
\Delta_n(\mathbf{x})\,
\left[ \prod_{i=1}^n\xi_i^{g-1}\right]\, p\,(\mathbf{x}),
\label{qop7}
\end{equation}
using the notation of Theorem \ref{dirichlet1}, notably
\begin{equation}
\xi_i=\frac{\prod_{k=1}^n(x_k-y_i)}{\prod_{k\neq i}(y_k-y_i)}\,.
\label{liov5bis}
\end{equation}

By virtue of the modification \eqref{liov13}
of the Dirichlet-Liouville integral, $Q_z$ is normalised
such that it sends $p\,(\mathbf{x})\equiv 1$ to $1$:
\begin{equation}
Q_z[1]=1.
\label{qop8}
\end{equation}
The function $[Q_zp](\mathbf{y})$ defined by the above integral
 for $\forall p\in\Cbbd[\mathbf{x}]^{S_n}$
as a real-valued function of real arguments $z$ and $\mathbf{y}$
is, in fact, a polynomial in $z$ and a symmetric polynomial
in $\mathbf{y}$.

\begin{prop}\label{prop2}
The operator $Q_z$ defined by \eqref{qop7} sends $\Cbbd[\mathbf{x}]^{S_n}$
into $\Cbbd(g)[\mathbf{y},z]^{S_n^y}$.
\end{prop}

{\bf Proof.}
Let $\mbox{\bf k}=(k_1,\dots,k_n)\in \Nbbd^n$.
Extend the domain of $Q_z$ from $\Cbbd[\mathbf{x}]^{S_n}$ to
$\Cbbd(\mathbf{y})[\mathbf{x}]^{S_n}$, assuming that $Q_z$ acts trivially on
the $\mathbf{y}$-variables. In other words, we can allow $p$ in the formula
\eqref{qop7} to depend on variables $\mathbf{y}$ in addition to $\mathbf{x}$.
In particular, we can set $p={\boldsymbol{\xi}}^{\mathbf{k}}$.
The integral \eqref{liov13} is well suited for calculating the action
of $Q_z$ on the monomials
${\boldsymbol{\xi}}^{\mathbf{k}}\equiv
\xi_1^{k_1}\cdots\xi_n^{k_n}\in\Cbbd(\mathbf{y})[\mathbf{x}]^{S_n}$.
Substituting $p={\boldsymbol{\xi}}^{\mathbf{k}}$ into \eqref{qop7}
and using
the integral \eqref{liov13} for $\a_i=g+k_i$,
we arrive at the following formula:
\begin{equation}
Q_z:\;{\boldsymbol{\xi}}^{\mathbf{k}}
\mapsto
{\boldsymbol{y}}^{\mathbf{k}}\,
\frac{(z-1)^{\abs{\mathbf{k}}}}{(ng)_{\abs{\mathbf{k}}}}\,
\prod_{i=1}^n(g)_{k_i},
\label{qop11}
\end{equation}
where
\begin{equation}
(\a)_k\equiv\frac{\Gamma(\a+k)}{\Gamma(\a)}
=\a(\a+1)\cdots(\a+k-1)\qquad \forall k\in \Nbbd,
\label{qop11a}
\end{equation}
is the Pochhammer symbol.

The monomials ${\boldsymbol{\xi}}^{\mathbf{k}}$ form a basis
in $\Cbbd[\mathbf{x}]^{S_n}$.
To show this, it is sufficient to express the elementary symmetric functions
$e_i(\mathbf{x})$, $i=1,\dots,n$ (cf \eqref{qop12})
in terms of $\boldsymbol{\xi}$, for
any symmetric polynomial in ${\bf x}$ is uniquely expanded into
the powers
$E_{\l}(\mathbf{x})
  =e_1^{\l_1-\l_2}(\mathbf{x})e_2^{\l_2-\l_3}(\mathbf{x})
\cdots e_n^{\l_n}(\mathbf{x})$ (see section \ref{jack}).
Introducing the generating function $w_x(t)$ for $\{e_i(\mathbf{x})\}$,
\begin{equation}
 w_x(t)=\prod_{j=1}^n(t-x_j)=t^n+
 \sum_{i=1}^n (-1)^i t^{n-i} e_i(\mathbf{x}),
\label{defw}
\end{equation}
and using the formula \eqref{Rres} for the rational function
$R(t)=w_x(t)/w_y(t)$, we get
\begin{equation}
 w_x(t)=w_y(t)R(t)
   =w_y(t)-\sum_{i=1}^n\xi_i\,\frac{w_y(t)}{t-y_i}\,.
\label{qop13}
\end{equation}
Expanding \eqref{qop13} in powers of $t$, we obtain the expressions
for $\{e_i(\mathbf{x})\}$ in terms of $\boldsymbol{\xi}$:
\begin{equation}
e_i(\mathbf{x})=e_i(\mathbf{y})
+\sum_{j=1}^n \xi_j \,\left[e_{i-1}(\mathbf{y})\right]_{y_j=0}.
\label{qop14}
\end{equation}

We can give the following algorithm for calculating
$[Q_zp](\mathbf{y})$ for an arbitrary symmetric
polynomial $p\in\Cbbd[\mathbf{x}]^{S_n}$:
\begin{enumerate}
\item express $p(\mathbf{x})$ as a polynomial in $\{e_i(\mathbf{x})\}$,
\item substitute the expressions \eqref{qop14} for $\{e_i(\mathbf{x})\}$ in terms
of $\{e_j(\mathbf{y})\}$ and $\boldsymbol{\xi}$,
\item expand the resulting polynomial in monomials
${\boldsymbol{\xi}}^{\mathbf{k}}$,
\item replace each monomial ${\boldsymbol{\xi}}^{\mathbf{k}}$
with the right hand side of \eqref{qop11}.
\end{enumerate}

The result $[Q_zp](\mathbf{y})$ 
is obviously a polynomial in $\mathbf{y}$ and $z$
and a rational function in $g$. Any permutation of the
$\mathbf{y}$-variables,
$\sigma:y_i\mapsto y_{\sigma_i}$, by \eqref{liov5},
also sends $\xi_i$ to $\xi_{\sigma_i}$. 
The right hand side of
the \eqref{qop14} is therefore invariant
under simultaneous permutations 
$\sigma:y_i\mapsto y_{\sigma_i}$ and $\xi_i\mapsto\xi_{\sigma_i}$.
Using \eqref{qop11}, we conclude that the polynomial $[Q_zp](\mathbf{y})$
is symmetric in $\mathbf{y}$.
\endproof

\begin{remark}\label{algorithm}
Let us mention another convenient variant 
of the above algorithm.
After the substitution $\xi_j=\eta_j y_j$, the expression \eqref{qop14}
for $e_i(\mathbf{x})$ becomes homogeneous in $\mathbf{y}$:
\begin{equation}
e_i(\mathbf{x})=\sum_{1\leq k_1<\ldots<k_i\leq n}
(1+\eta_{k_1}+\ldots+\eta_{k_i})y_{k_1}\cdots y_{k_i}. 
\label{qop36}
\end{equation}

We can use \eqref{qop36} instead of \eqref{qop14} for the step 2 of the
algorithm, and then use the substitution
\begin{equation}
B_g: \;\boldsymbol{\eta}^{\mathbf{k}}
\mapsto
\frac{(z-1)^{\abs{\mathbf{k}}}}{(ng)_{\abs{\mathbf{k}}}}
\prod_{i=1}^n(g)_{k_i},
\label{subs-eta}
\end{equation}
instead of \eqref{qop11}, for the step 4.
\end{remark}

\begin{remark}\label{algorithm-x}
Notice also that the polynomiality of $Q_z\boldsymbol{\xi^k}$ in $z$ 
\eqref{qop11} ensures that
the operator $Q_z$, although initially defined for the real values $z>1$, 
can now be analytically continued to the whole of the complex plane,
$z\in\Cbbd$.
\end{remark}

Formula \eqref{qop7} is not the only way to realise $Q_z$ as an integral
operator. An alternative description of $Q_z$ can be obtained
{} from Theorem \ref{dirichlet2}.

For $\gamma\in(0,\pi)$ and for an ordered set
of real variables $\mathbf{v}=(v_1,\dots,v_n)$ satisfying
\begin{displaymath}
 v_1< v_2<\ldots< v_n< \pi+v_1,
\end{displaymath}
we define the operator $\widetilde{Q}_\gamma$ acting on
a symmetric polynomial $p\in\Cbbd[\mathbf{x}]^{S_n}$
by the integral
\begin{multline}
 [\widetilde{Q}_\gamma p](\mathbf{v})=
\frac{\Gamma(ng)}{\Gamma(g)^n}
\,\frac{(\sin\gamma)^{1-ng}}{\ds\prod_{i<j}\sin(v_j-v_i)} \\
\times\int_{\Omega^T_{\mathbf{u}}}d\mathbf{u}\;
\delta(\abs{\mathbf{u}}-\abs{\mathbf{v}}-\gamma)
\left[\prod_{j<k}\sin(u_k-u_j)\right]
\left[\prod_{j=1}^n(\xi^T_j)^{g-1}\right]\,
p\,(\varepsilon(\mathbf{u})),
\label{qop17}
\end{multline}
using the notation 
of Theorem \ref{dirichlet2}
and the map $\varepsilon$ from \eqref{jack18a}.

\begin{prop}
Let
\begin{equation}
 \mathbf{x}=\varepsilon(\mathbf{u}),
\qquad 
\mathbf{y}=\varepsilon(\mathbf{v}),
\qquad z=\varepsilon(\gamma).
\label{uvg-xyz}
\end{equation}
Then $[\widetilde{Q}_\gamma p](\mathbf{y})\in\Cbbd(g)[\mathbf{y},z]^{S_n^y}$
and, moreover, $\widetilde{Q}_\gamma$ coincides with $Q_z$ as operators from
$\Cbbd[\mathbf{x}]^{S_n}$ to $\Cbbd[\mathbf{y}]^{S_n}$.
\label{equivQ}
\end{prop}

{\bf Proof.}
Substituting $p\left(\varepsilon(\mathbf{u})\right)=\left(\boldsymbol{\xi}^T\right)^{\mathbf{k}}$ into
\eqref{qop17} and using the integral \eqref{liov23}, one obtains the action
\begin{equation}
\widetilde Q_\gamma:
\left(\boldsymbol{\xi}^T\right)^{\mathbf{k}}
\mapsto
\frac{(\sin\gamma)^{\abs{\mathbf{k}}}}{(ng)_{\abs{\mathbf{k}}}}
\prod_{i=1}^n(g)_{k_i},
\label{qop11b}
\end{equation}
which is similar to \eqref{qop11}. Comparing expressions \eqref{liov5} for $\xi_i$
and \eqref{liov17} for $\xi_i^T$ and using substitutions \eqref{uvg-xyz},
we obtain the relation
\begin{equation}
  \xi_j=\xi_j^T\cdot 2iy_je^{i\gamma},
\label{relxi}
\end{equation}
which identifies the action \eqref{qop11} of $Q_z$ on the monomials
$\boldsymbol{\xi}^{\mathbf{k}}$ and the action \eqref{qop11b} of $\widetilde Q_\gamma$
on the monomials $\left(\boldsymbol{\xi}^T\right)^{\mathbf{k}}$,
and therefore on polynomials $p\in\Cbbd[\mathbf{x}]^{S_n}$.
Hence, $Q_z$ and $\widetilde Q_\gamma$ 
define two integral realisations of 
the same operator acting from $\Cbbd[\mathbf{x}]^{S_n}$ 
to $\Cbbd[\mathbf{y}]^{S_n}$.
\endproof

The integral formula \eqref{qop7} will be our main tool when working with
the $Q$-operator. The trigonometric representation \eqref{qop17}
 will, however, be required for proving the conjugation properties of 
$Q_z$ (cf Proposition \ref{prop3} below). In what follows, we shall 
identify the source space $\Cbbd[\mathbf{x}]^{S_n}$
and the target space $\Cbbd[\mathbf{y}]^{S_n}$ of $Q_z$ and think of
$Q_z$ as an operator in $\Cbbd[\mathbf{x}]^{S_n}$.

\begin{prop}\label{commQe}
Let $\hat{e}_n$ be the multiplication operator acting
in $\Cbbd[\mathbf{x}]^{S_n}$ as
$\hat{e}_n:p\,(\mathbf{x})\mapsto e_n(\mathbf{x})\,p\,(\mathbf{x})$.
Then the following operator relation holds:
\begin{equation}
  Q_z\hat{e}_n=z\hat{e}_nQ_z.
\label{Qe=zeQ}
\end{equation}
\end{prop}
{\bf Proof.} For $z>1$, the statement follows from the observation
that $Q_z$ is realised as the integral
operator \eqref{qop7} whose kernel contains
$\delta\left(e_n(\mathbf{x})-ze_n(\mathbf{y})\right)$. Then the result is extended by
analyticity to $z\in\Cbbd$.
\endproof

The following Theorem is the main result of this section, because it states
the defining properties \eqref{qop2}--\eqref{qop4} of the introduced 
$Q$-operator.

\begin{theo}\label{qop41}
Jack polynomials are the eigenfunctions of the operator $Q_z$:
\begin{equation}
[Q_z P_\l^{(1/g)}](\mathbf{x})=q_\l(z)P_\l^{(1/g)}(\mathbf{x}),\qquad
  q_\l(z)\in\Cbbd(g)[z].
\label{qop42}
\end{equation}
\end{theo}

The idea of the proof is to show that the matrix of $Q_z$ in the basis
of Jack polynomials is both symmetric
and triangular, and therefore diagonal. The symmetry is established by
Proposition \ref{prop3}, the triangularity by Proposition \ref{qop38}.

\begin{prop}\label{prop3}
The operator $Q_z$ possesses the conjugation property
\begin{equation}
  Q_z^\dagger=Q_{\bar{z}}
\label{conjugQ}
\end{equation}
with respect to the Hermitian form \eqref{jack20}.
\end{prop}

{\bf Proof.}
We shall use the representation of the $Q$-operator as $\widetilde{Q}_\gamma$
given by \eqref{qop17}. To begin with, let us prove that
\begin{equation}
 \langle f,\widetilde Q_\gamma p\,\rangle=\langle \widetilde Q_{\pi-\gamma}f,p\,\rangle\qquad
 \forall f,p\in\Cbbd[\mathbf{x}]^{S_n}, \qquad \gamma\in(0,\pi).
\label{qop21}
\end{equation}

For our present purposes it is convenient to think that,
similarly to \eqref{defTv}, the variables $\mathbf{u}$, $\mathbf{v}$
in \eqref{qop17} belong to the torus $T_n=\Rbbd^n/\pi\Zbbd^n$
and all functions of $\mathbf{u}$, $\mathbf{v}$, including
delta-function, are $\pi$-periodic.
In particular, we can assume that all the sine functions in \eqref{qop17}
take their arguments in $(0,\pi)$ or, equivalently, replace the sines
by their absolute values.

Keeping that in mind and using \eqref{jack20} and \eqref{qop17}
we obtain
\begin{equation}
  \langle f,\widetilde Q_\gamma p\,\rangle=\int_{T_{\mathbf{u},\mathbf{v}}} d\mathbf{u}\,d\mathbf{v}\,
  \overline{f(\varepsilon(\mathbf{v}))}\,
  \Qcal_\gamma(\mathbf{v},\mathbf{u})
 \, p\,(\varepsilon(\mathbf{u})),
\label{qop24}
\end{equation}
where
\begin{equation}
  T_{\mathbf{u},\mathbf{v}}=\{\mathbf{u},\mathbf{v}\in\Rbbd^n/\pi\Zbbd^n \mid
   u_j\in(v_j,v_{j+1}), \quad v_j\in(u_{j-1},u_j)\}
\label{defTuv}
\end{equation}
(using the convention $j+n\equiv j$) and
\begin{eqnarray}
  \Qcal_\gamma(\mathbf{v},\mathbf{u}) &=&
  (\sin\gamma)^{1-ng}\,\frac{\Gamma(ng)}{\Gamma(g)^n}\;
  \delta(|\mathbf{u}|-|\mathbf{v}|-\gamma) \nonumber \\
  &&\times \prod_{j<k}\left|\sin(v_k-v_j)\sin(u_k-u_j)\right|
  \prod_{j,k}\left|\sin(u_k-v_j)\right|^{g-1}.
\label{Quv}
\end{eqnarray}

Performing the following change of variables:
\begin{equation}
  u_j\mapsto u_j'=v_{j+1},\qquad v_j\mapsto v_j'=u_{j},
\label{qop26}
\end{equation}
and observing that
$T_{\mathbf{u},\mathbf{v}}=T_{\mathbf{u}',\mathbf{v}'}$ and
$\Qcal_\gamma(\mathbf{v},\mathbf{u})
=\Qcal_{\pi-\gamma}(\mathbf{v}',\mathbf{u}')$
we transform \eqref{qop24} into $(\widetilde Q_{\pi-\gamma}f,p)$.

Having thus proved that $\widetilde{Q}_\gamma^\dagger=\widetilde{Q}_{\pi-\gamma}$
for $\gamma\in(0,\pi)$, we extend this result by analyticity to
$\widetilde{Q}_\gamma^\dagger=\widetilde{Q}_{\pi-\bar\gamma}$ for $\gamma\in\Cbbd$.
Due to the relation $z=e^{2i\gamma}$ (cf \eqref{uvg-xyz}), the involution
$\gamma\mapsto\pi-\bar\gamma$ corresponds to $z\mapsto\bar{z}$,
hence \eqref{conjugQ}.
\endproof

\begin{prop}\label{qop38}
The matrix of the operator $Q_z$ between any two of the three bases,
$m_\lambda$ \eqref{def-m}, $E_\lambda$ \eqref{qop37} or $P_\lambda^{(1/g)}$,
is triangular with respect to the dominance partial ordering.
For example,
\begin{equation}
Q_z E_\l=\sum_{\mu\preceq\l}q_{\l\mu}(z,g)m_\mu.
\label{qop39}
\end{equation}
\end{prop}

{\bf Proof.}
Since the three bases, $\{m_\lambda, E_\lambda, P_\lambda^{(1/g)}\}$,
are mutually triangular, it is sufficient to prove only the expansion
\eqref{qop39}.

According to the algorithm explained in Remark \ref{algorithm}, in order
to calculate $Q_z E_{\lambda}$ we have
to substitute into \eqref{qop37} the expressions \eqref{qop36} for $e_j$
in terms of $\eta_j$ and $y_j$, expand the brackets
and transform the products of $\eta_j$ by the formula \eqref{subs-eta}.

It follows from \eqref{qop36} that $Q_z$ maps
$E_{\l}$ into a polynomial consisting of the same monomials
in $y_j$ as $E_{\l}$ in $x_j$. The only difference is the presence of
the variables $\eta_j$ in the coefficients producing a nontrivial dependence
on $g$. Since the expansion of $E_{\lambda}$ into $m_\mu$ is triangular
\eqref{qop40}, so is the expansion of $Q_zE_{\lambda}$.
\endproof

{\bf Proof of Theorem \ref{qop41}.}
Let $\Qcal_{\l\mu}\equiv(Q_zP^{(1/g)}_\l,P^{(1/g)}_\mu)$.
It follows from the triangularity property (Proposition \ref{qop38})  that
$\Qcal_{\l\mu}\neq0$ only for $\l\succeq\mu$. On the other hand,
{} from the conjugation formula \eqref{conjugQ} and the orthogonality of
the basis $P^{(1/g)}_\mu$ it follows that
$(\Qcal_{\l\mu}\neq0)\Longrightarrow(\Qcal_{\mu\l}\neq0)$, therefore
$\mu\succeq\l$. However, $\l\succeq\mu$ and $\mu\succeq\l$ implies
$\l=\mu$, hence \eqref{qop42}.
The eigenvalue $q_\lambda(z)$ is a polynomial in $z$ since
$[Q_zP^{(1/g)}_\l](\mathbf{x})\in\Cbbd(g)[x,z]^{S_n^x}$ by Proposition \ref{prop2}.
\endproof

\begin{remark} In the quasi-classical limit,
$g\rightarrow\infty$, the operator $Q_z$ turns into the
corresponding B\"acklund transformation for the
Calogero-Sutherland model (see \cite{KS5} for details).
\end{remark}

\section{Eigenvalues of $Q_z$}\label{eigenvalues}
\noindent
Here we study
the eigenvalues of the operator $Q_z$,
that is the polynomials $q_\lambda(z)\in\Cbbd(g)[z]$
introduced in the previous section (cf  \eqref{qop42}).

First, we derive a representation for $q_\lambda(z)$
as a finite multiple sum by making use of the
explicit formula for the action of $Q_z$ on a polynomial basis.
Then, we convert the obtained expression into
another finite multiple sum, which results from a
previously known formula (cf \cite{KS4,KS96}) for the
separated polynomials $q_\lambda(z)$ in terms
of the hypergeometric function ${}_nF{}_{n-1}(z)$.
Finally, we write down the ordinary
differential equation of order $n$ for $q_\lambda(z)$,
which serves as the Baxter equation \eqref{qop5}
for our $Q$-operator $Q_z$.

By Proposition \ref{qop38}, the eigenvalue $q_\l(z)$
can be determined
as the diagonal matrix element of $Q_z$ in the basis $E_{\l}$:
\begin{equation}
  Q_zE_{\l}=q_\l(z)E_{\l}
  +\sum_{\mu\prec\l}b_{\l\mu}E_{\mu}\,.
\label{q-diag}
\end{equation}
To calculate $Q_zE_{\l}$ we shall apply the algorithm described
in Remark \ref{algorithm}.
Taking the diagonal element corresponds to
replacing $e_j(\mathbf{x})$ by formula \eqref{qop36} and
taking the coefficient at the monomial $\mathbf{y}^{\l}$.
Introducing the operator $B_g:\Cbbd[\boldsymbol{\eta}]\mapsto\Cbbd(g)$
acting on the monomials $\boldsymbol{\eta}^{\mathbf{k}}$ by the formula
\eqref{subs-eta} we arrive at the following expression:
\begin{equation}
q_\l(z)=B_g
\left[\prod_{i=1}^n(1+\eta_1+\ldots+\eta_i)^{\lambda_{i,i+1}}\right].
\label{qop49}
\end{equation}
\begin{theo}\label{theorem4}
The polynomial $q_\l(z)\in\Cbbd(g)[z]$ is given by the multiple sum
\begin{equation}
q_\l(z)=
 z^{\l_n}\sum_{k_1=0}^{\lambda_{1,2}}\cdots\sum_{k_{n-1}=0}^{\lambda_{n-1,n}}
\prod_{i=1}^{n-1}(1-z)^{k_i}\,\frac{(-\lambda_{i,i+1})_{k_i}}{k_i!}\,
\frac{(ig)_{k_1+\ldots+k_i}}{((i+1)g)_{k_1+\ldots+k_i}}\,.
\label{qop47}
\end{equation}
\end{theo}

{\bf Proof.}
Expanding the product in (\ref{qop49}) we obtain
\begin{equation}
q_\l(z)=\sum_{k_1=0}^{\lambda_{1,2}}\cdots\sum_{k_n=0}^{\lambda_{n,n+1}}
\left[\,\prod_{j=1}^n\binom{\lambda_{j,j+1}}{k_j}\right]
\frac{(z-1)^{\abs{\mathbf{k}}}}{(ng)_{\abs{\mathbf{k}}}}\;\,
c_{\mathbf{k}}
\label{exprq1}
\end{equation}
with
\begin{eqnarray}
c_{\mathbf{k}}&=&\frac{(ng)_{\abs{\mathbf{k}}}}
{(z-1)^{\abs{\mathbf{k}}}}\;B_g
\left[\prod_{j=1}^n
(\eta_1+\ldots+\eta_j)^{k_j}\right] \nonumber\\
&=& \left(D_1\otimes\cdots\otimes D_n\right)
\left[\prod_{j=1}^n
(\eta_1+\ldots+\eta_j)^{k_j}\right],
\label{qop48}
\end{eqnarray}
where the operators $D_j:\Cbbd[\eta_j]\mapsto\Cbbd(g)$ are defined
as follows $D_j:\eta_j^{k_j}\mapsto(g)_{k_j}$.

\begin{lemma}\label{lemmaq1}
The coefficients $c_{\mathbf{k}}$ \eqref{qop48} are
\begin{equation}
c_{\mathbf{k}}=\prod_{j=1}^n (jg+k_1+\ldots+k_{j-1})_{k_j}=
(ng)_{|\mathbf{k}|}\prod_{j=1}^{n-1}
\frac{(jg)_{k_1+\ldots+k_j}}{((j+1)g)_{k_1+\ldots+k_j}}\,.
\label{qcoeff2}
\end{equation}
\end{lemma}
The formula \eqref{qop47} is then obtained by
substituting \eqref{qcoeff2} into \eqref{exprq1} and evaluating the sum
over $k_n$.
\endproof

{\bf Proof of Lemma \ref{lemmaq1}} is given by induction in $n$.
When $n=1$ formulae \eqref{qop48} and \eqref{qcoeff2} give the same value
$(g)_{k_1}$ for $c_{k_1}$.
It is sufficient to show that
$c_{k_1,\dots,k_n}$ given by \eqref{qop48} satisfies the recurrence relation
\begin{equation}
c_{k_1,\dots,k_{n}}=c_{k_1,\dots,k_{n-1}}
\cdot  \left(ng+k_1+\ldots+k_{n-1}\right)_{k_{n}}.
\label{recurr-a}
\end{equation}
Expanding the factor $((\eta_1+\ldots+\eta_{n-1})+\eta_n)^{k_n}$ in
\eqref{qop48} we get
\begin{eqnarray}
c_{k_1,\dots,k_{n}}&=&
\left(D_1\otimes\cdots\otimes D_{n}\right)
\left[\biggl(\,\prod_{j=1}^{n-1}(\eta_1+\ldots+\eta_j)^{k_j}\biggr)
(\eta_1+\ldots+\eta_n)^{k_n} \right]\nonumber \\
&=& \sum_{m=0}^{k_n}\binom{k_n}{m}
\left(D_1\otimes\cdots\otimes D_{n-1}\right)\biggl[
(\eta_1+\ldots+\eta_{n-1})^{k_{n-1}+m}\biggr. \nonumber \\
&&\times\biggl.\prod_{j=1}^{n-2}(\eta_1+\ldots+\eta_j)^{k_j}\biggr]\;\cdot
D_n\left[\eta_n^{k_{n}-m}\right] \nonumber \\
&=& \sum_{m=0}^{k_n}\binom{k_n}{m}\;
c_{k_1,\dots,k_{n-2},k_{n-1}+m} \;(g)_{k_n-m} \nonumber \\
&=& c_{k_1,\dots,k_{n-1}}\,
\sum_{m=0}^{k_n}\binom{k_n}{m}
((n-1)g+k_1+\ldots+k_{n-1})_{m}(g)_{k_n-m}.\qquad
\end{eqnarray}
To obtain \eqref{recurr-a} it remains to use the summation formula
\begin{equation}
  \sum_{m=0}^k \binom{k}{m}(s)_m(g)_{k-m}=(s+g)_k
\label{ident}
\end{equation}
for $k=k_n$ and $s=(n-1)g+k_1+\ldots+k_{n-1}$ which
is easily proved
by comparing the coefficients with powers of $t$ in the identity 
$(1-t)^{-s}(1-t)^{-g}=(1-t)^{-s-g}$.
\endproof

In \cite{KS4} (see also \cite{KS96}) it was conjectured,
and then proved for $n=2,3,4$ (cf \cite{KS4,KS96,Man,Man2}), 
that Jack polynomials $P_\l^{(1/g)}(x_1,\dots,x_n)$
can be factorised into a product of the separated polynomials 
$f_\l(z)$ defined in terms of the hypergeometric function
${}_nF{}_{n-1}(z)$ as follows:
\begin{equation}
f_\l(z)=z^{\l_n}(1-z)^{1-ng}\phantom{|}_nF_{n-1}\left(
\begin{matrix} a_1,a_2,\dots,a_n\cr
b_1,\dots,b_{n-1}
\end{matrix};z
\right),
\label{conn2}
\end{equation}
\begin{equation}
a_i=\l_{n,i}+1-(n-i+1)g,\quad b_i=a_i+g. \label{conn4}
\end{equation}
\begin{theo}\label{theorem1}
The function $f_\l(z)$ is a polynomial given by the multiple sum
\begin{eqnarray}
f_\l(z)&=&z^{\l_n}\sum_{k_1=0}^{\lambda_{1,2}}\cdots\sum_{k_{n-1}=0}^{\lambda_{n-1,n}}
\prod_{i=1}^{n-1}(1-z)^{\lambda_{i,i+1}}\Bigl(\frac{z}{z-1}\Bigr)^{k_i}
\frac{(-\lambda_{i,i+1})_{k_i}}{k_i!} \nonumber\\
&&\times
\frac{(1-(n-i+1)g-\lambda_{i,n})_{k_i+\ldots+k_{n-1}}}%
{(1-(n-i)g-\lambda_{i,n})_{k_i+\ldots+k_{n-1}}}\,.
\label{conn5}
\end{eqnarray}
\end{theo}
{\bf Proof.} Representation \eqref{conn5} for 
the separated polynomial $f_\l(z)$ follows from 
the observation made in \cite{KS96} that an
explicit {\it finite} multiple sum for the polynomial 
given by \eqref{conn2}
can be obtained by iterating Fox's formula \cite{Fox},
\begin{multline}
\phantom{|}_pF_{p-1}\left(
\begin{matrix}
a_1,\dots,a_{p-1},b_{p-1}+m\cr
b_1,\dots,b_{p-1}
\end{matrix}
;z\right)
=\sum_{j=0}^m(-z)^j\,\frac{(-m)_j}{j!} \\
\times \left[\prod_{i=1}^{p-1}
\frac{(a_i)_j}{(b_i)_j}\right]\phantom{|}_{p-1}F_{p-2}\left(
\begin{matrix}
a_1+j,\dots,a_{p-1}+j\cr
b_1+j,\dots,b_{p-2}+j
\end{matrix}
;z\right),\qquad
\label{conn6}
\end{multline}
which expresses the hypergeometric function
$\phantom{|}_pF_{p-1}$ on the left, for any 
non-negative integer $m$,
as a {\it finite} sum of lower series.

Noticing that $a_{i+1}=b_i+\lambda_{i,i+1}$ for $i=1,\dots,n-1$ in
\eqref{conn4} and using the binomial theorem, ${}_1F_0\left(
\begin{matrix}
a\cr-
\end{matrix}
;z\right) = (1-z)^{-a}$,
we apply \eqref{conn6} to \eqref{conn2} consecutively for $p=n-1,\dots,2$
and derive  \eqref{conn5}.
\endproof
\begin{theo}\label{theorem2}
The polynomials $f_\l(z)$ and $q_\l(z)$ are proportional
\begin{equation}
f_\l(z)=\beta_\lambda q_\l(z),
\label{conn18}
\end{equation}
where
\begin{equation}
\beta_\lambda=\prod_{i=1}^{n-1}\beta_\lambda^{(i)}, \qquad
\beta_\lambda^{(i)}=\frac{((n-i+1)g)_{\l_{i,n}}}{((n-i)g)_{\l_{i,n}}}\,.
\label{conn19}
\end{equation}
\end{theo}

{\bf Proof.} We shall start with the expression \eqref{conn5} for $f_\l(z)$.
After the change of the variable $z\mapsto1-z$ we can rewrite \eqref{conn5}
as
\begin{equation}
  f_\l(1-z)=(1-z)^{\l_n}\varphi_\lambda^{(n-1)}(z),
\label{conn20}
\end{equation}
where $\varphi_\lambda^{(i)}(z;\allowbreak k_{i+1},\dots,k_{n-1})$
are defined recursively for $i=1,\dots,n-1$ by
$\varphi_\lambda^{(0)}(z;\allowbreak k_1,\allowbreak\dots,k_{n-1})\equiv1$ and
\begin{multline}
\varphi_\lambda^{(i)}(z;k_{i+1},\dots,k_{n-1})
=z^{\lambda_{i,i+1}}\sum_{k_i=0}^{\lambda_{i,i+1}}
\left(\frac{z-1}{z}\right)^{k_i}
\frac{(-\lambda_{i,i+1})_{k_i}}{k_i!} \\
\times\frac{(1-(n-i+1)g-\lambda_{i,n})_{k_i+\ldots+k_{n-1}}}%
{(1-(n-i)g-\lambda_{i,n})_{k_i+\ldots+k_{n-1}}}\;
\varphi_\lambda^{(i-1)}(z;k_{i},\dots,k_{n-1}).\qquad
\label{conn20a}
\end{multline}
\begin{lemma}
\begin{multline}
\varphi_\lambda^{(i)}(z;k_{i+1},\dots,k_{n-1})=
\beta_\l^{(i)}\,\frac{\Gamma((n-i+1)g)}{\Gamma(g)\Gamma((n-i)g)} \\
\times \int_0^1 dv\, v^{ig-1}
(1-v)^{(n-i)g-1+\lambda_{i+1,n}-k_{i+1}-\ldots-k_{n-1}}
(1-zv)^{\lambda_{i,i+1}}
\Phi_\lambda^{(i-1)}(zv),
\label{phiPhi}
\end{multline}
where $\Phi_\lambda^{(i)}(z)$  for $i=0,\dots,n-2$
are defined recursively by  $\Phi_\lambda^{(0)}(z)\equiv1$ and
\begin{equation}
 \Phi_\lambda^{(i)}(z)
=\int_0^1 du\, u^{ig-1}(1-u)^{g-1}(1-zu)^{\lambda_{i,i+1}}
\Phi_\lambda^{(i-1)}(zu).
\label{defPhi}
\end{equation}
\end{lemma}

{\bf Proof} is done  by induction. For $i=1$, formula \eqref{conn20a} gives
\begin{multline}
\varphi_\lambda^{(1)}(z;k_{2},\dots,k_{n-1})=
z^{\l_{1,2}}\,\frac{(1-ng-\l_{1,n})_{k_2+\ldots+k_{n-1}}}
{(1-(n-1)g-\l_{1,n})_{k_2+\ldots+k_{n-1}}} \\
\times\,{}_2F{}_1\left(
\begin{matrix}
-\l_{1,2},\theta_1\cr
\theta_1+g
\end{matrix};\frac{z-1}{z}
\right),\qquad\qquad\qquad
\end{multline}
where $\theta_1=1-ng-\l_{1,n}+k_2+\ldots+k_{n-1}$.
Using the well-known identity for the terminating Gauss
hypergeometric function,
\begin{equation}
\phantom{|}_2F_1\left(
\begin{matrix}
-p,b\cr c
\end{matrix}
;z\right)
=\frac{(b)_p}{(c)_p}\;(1-z)^p\phantom{|}_2F_1\left(
\begin{matrix}
-p,c-b
\cr -p+1-b
\end{matrix}
;\frac{1}{1-z}\right),
\label{conn8}
\end{equation}
and then the integral representation (see, e.g. \cite{Bateman}),
\begin{equation}
\phantom{|}_2F_1\left(
\begin{matrix}
a,b\cr c
\end{matrix}
;z\right)=\frac{\Gamma(c)}{\Gamma(b)\Gamma(c-b)}\int_0^1du\>u^{b-1}(1-u)^{c-b-1}
(1-zu)^{-a},
\label{iii}
\end{equation}
we arrive at the formula \eqref{phiPhi} for $i=1$.

Supposing that the statement is true for $i-1$, we substitute the
expression \eqref{phiPhi} for $\varphi_\lambda^{(i-1)}$ in terms of $\Phi_\lambda^{(i-2)}$
 into \eqref{conn20a} and get the integral
\begin{multline}
 \varphi_\lambda^{(i)}(z;k_{i+1},\dots,k_{n-1})
= z^{\lambda_{i,i+1}} \\
 \times\int_0^1 dv\, v^{ig-1}
(1-v)^{(n-i+1)g-1+\lambda_{i,n}-k_{i+1}-\ldots-k_{n-1}}
(1-zv)^{\lambda_{i-1,i}}\Phi_\lambda^{(i-2)}(zv)S_i,
\label{indstep1}
\end{multline}
where $S_i$ is the sum
\begin{multline}
S_i=\sum_{k_i=0}^{\lambda_{i,i+1}}
\left(\frac{z-1}{z(1-v)}\right)^{k_i}\,
\frac{(-\lambda_{i,i+1})_{k_i}}{k_i!}\,
\frac{(1-(n-i+1)g-\lambda_{i,n})_{k_i+\ldots+k_{n-1}}}%
{(1-(n-i)g-\lambda_{i,n})_{k_i+\ldots+k_{n-1}}} \\
=\frac{(1-(n-i+1)g-\lambda_{i,n})_{k_{i+1}+\ldots+k_{n-1}}}%
{(1-(n-i)g-\lambda_{i,n})_{k_{i+1}+\ldots+k_{n-1}}}
\phantom{|}_2F_1\left(
\begin{matrix}
-\lambda_{i,i+1},\theta_i\cr\theta_i+g
\end{matrix}
;\frac{(z-1)}{z(1-v)}\right)
\label{conn21}
\end{multline}
and
$\theta_i=1-(n-i+1)g-\lambda_{i,n}+k_{i+1}+\ldots+k_{n-1}$.
Using again the transformation \eqref{conn8} and then 
the integral representation \eqref{iii},
after some algebra we transform $S_i$ to the following form:
\begin{multline}
S_i=\beta_\lambda^{(i)}\,\frac{\Gamma((n-i+1)g)}{\Gamma(g)\Gamma((n-i)g)}\;
(1-v)^{-\lambda_{i,i+1}} \\
\times\int_0^1dt\,
t^{g-1}(1-t)^{(n-i)g-1+\lambda_{i+1,n}-k_{i+1}-\ldots-k_{n-1}}
\bigl[1-y(v+t(1-v))\bigr]^{\lambda_{i,i+1}}.
\label{Sformula}
\end{multline}
Substituting $S_i$ from \eqref{Sformula} into \eqref{indstep1} we obtain
an expression for $\varphi_\lambda^{(i)}$ as a double integral, in the variables $v$ and $t$. 
Making the change of variables
$(v,t)\mapsto(v,u=v+t-vt)$ followed by $(u,v)\mapsto(u,w=v/u)$
we obtain for $\varphi_\lambda^{(i)}$ an expression as an integral in
$u$ and $w$. The integration domains are changed respectively
{} from $v,t\in[0,1]$ to
$0\leq v\leq u\leq1$, and finally to $u,w\in[0,1]$.
The result reads
\begin{multline}
\varphi_\lambda^{(i)}(z;k_{i+1},\dots,k_{n-1})
 =\beta_\lambda^{(i)}\,\frac{\Gamma((n-i+1)g)}{\Gamma(g)\Gamma((n-i)g)} \\
\times \int_0^1 du\, u^{ig-1}
(1-u)^{(n-i)g-1+\lambda_{i+1,n}-k_{i+1}-\ldots-k_{n-1}}
(1-zu)^{\lambda_{i,i+1}} \\
\times \int_0^1 dw\, w^{(i-1)g-1}(1-w)^{g-1}(1-zuw)^{\lambda_{i-1,i}}
\Phi_\lambda^{(i-2)}(zuw).
\label{int-uw}
\end{multline}
The integral over $w$ is evaluated with the help of \eqref{defPhi}
producing $\Phi_\lambda^{(i-1)}(zu)$, and we recover formula \eqref{phiPhi},
thus completing the induction step.
\endproof

Iterating integrals \eqref{defPhi} and using \eqref{conn20} we obtain
the following representation for $f_\l(1-z)$ as a multiple integral:
\begin{multline}
f_\l(1-z)=(1-z)^{\l_n}\beta_\lambda\,\frac{\Gamma(ng)}{\Gamma(g)^n}\,
\int_{0}^1du_1\cdots\int_0^1du_{n-1} \\
\times\prod_{i=1}^{n-1}
u_i^{ig-1}(1-u_i)^{g-1}(1-zu_i\cdots u_{n-1})^{\lambda_{i,i+1}}.\qquad\qquad
\label{conn33}
\end{multline}
Expanding the binomials $(1-zu_i\cdots u_{n-1})^{\l_{i,i+1}}$
we get a sum of products of one-dimensional beta-integrals in each $u_i$ and,
after evaluating them, we recover the expression \eqref{qop47}
for $q_\l(z)$, arriving thus at \eqref{conn18}.
\endproof

\begin{remark} The eigenvalue $q_\l(z)$ of the $Q$-operator $Q_z$
is thus proportional to the function $f_\l(z)$ expressed by \eqref{conn2} as
a hypergeometric polynomial, hence its Baxter (or separation) 
equation has the form (cf, e.g., \cite{Bateman})
\begin{equation}
\left(z\frac{d}{dz}\prod_{i=1}^{n-1}\left(z\frac{d}{dz}+b_i-1\right)-
z\prod_{i=1}^n\left(z\frac{d}{dz}+a_i\right)\right)\left[
z^{-\lambda_n}(1-z)^{ng-1}f_\lambda(z)\right]=0,
\label{diffeq}
\end{equation}
where $a_i=\l_{n,i}+1-(n-i+1)g$ and $b_i=a_i+g$.
\end{remark}

\section{Factorisation of Jack polynomials}\label{restricted}
\noindent
In this section, the integral operator $Q_z$ is used to construct a
factorisation of Jack polynomials
$P_\l^{(1/g)}(x_1,\dots,x_n)$ in $n$ variables
into a product of $n$ polynomials in one variable.
The general construction of the separating (factorising) operator
$\mathcal{S}_n^{(\rho)}$,
\begin{equation}
\mathcal{S}_n^{(\rho)}=\rho Q_{z_1}\cdots Q_{z_n},
\label{defSn}
\end{equation}
briefly described in the Introduction
involves an arbitrary linear functional 
$\rho:\Cbbd[\mathbf{x}]^{S_n}\mapsto\Cbbd$. 
Our present task is to make
the integral operator $\mathcal{S}_n^{(\rho)}$ as simple as possible 
by choosing $\rho=\rho_0$
in a special way. 
More specifically, we wish to factorise
the integral operator $\mathcal{S}_n^{(\rho_0)}$ into a string
of operators $\mathcal{A}_k$, 
\begin{equation}
\mathcal{S}_n^{(\rho_0)}=\rho_0 Q_{z_1}\cdots Q_{z_n}
=\mathcal{A}_1\cdots\mathcal{A}_n,
\end{equation}
possessing the property of {\it dimension reduction}
in the sense that each operator $\Acal_k$ acts only on $k$ variables.
Any such factorisation will be called {\it factorised separation chain}
or {\it $\mathcal{A}$-chain}.

Let us choose $\rho_0$ to be the evaluation homomorphism that
returns the value of a polynomial at the point $\mathbf{x}=(1,\dots,1)$,
\begin{equation}
  \rho_0:f\mapsto \rho_0 f=f(1,\dots,1),
\label{defrho0}
\end{equation}
or, equivalently, the integral operator whose kernel is
the product of delta-functions,
\begin{equation}
   \rho_0(\mathbf{x})=\prod_{j=1}^n\delta(x_j-1), \qquad
    \rho_0 f \equiv
    \int d\mathbf{x}\, \rho_0(\mathbf{x})f(\mathbf{x})=
f(1,\dots,1).
\end{equation}
Having fixed  the functional $\rho=\rho_0$, hereafter we shall omit
$\rho$ (or $\rho_0$) from the operator notation, writing $\Scal_n$ to denote
the corresponding separating operator. Two theorems below
constitute the main result of this section and show how distinguished
is the chosen $\rho_0$ \eqref{defrho0}.

\begin{theo} \label{decompSx}
The operator $\Scal_n=\rho_0Q_{z_1}\cdots Q_{z_n}$
sends Jack polynomials into the product of
separated polynomials $q_\lambda(z)$ defined in Section \ref{eigenvalues},
\begin{equation}
   \Scal_n:P_\lambda^{(1/g)}(\mathbf{x})
   \mapsto c_\lambda\prod_{k=1}^n q_\lambda(z_k),
\label{SP=qqqP}
\end{equation}
the normalisation constant $c_\lambda$ being given by
\begin{equation}
    c_\lambda=P_\lambda^{(1/g)}(1,\dots,1)=\prod_{1\le i<j\le n}\frac{(g(j-i+1))_{\l_{i,j}}}
{(g(j-i))_{\l_{i,j}}}\,.
\label{def-clambda}
\end{equation}
\end{theo}

{\bf Proof.}
Apply the chain \eqref{SP=qqqP} of operators $Q_{z_k}$ to a Jack
polynomial $P_\lambda^{(1/g)}(\mathbf{x})$ and use the fact \eqref{qop42}
that $P_\lambda^{(1/g)}(\mathbf{x})$ is an eigenfunction of $Q_{z_k}$.
Each factor $Q_{z_k}$ produces a desired factor $q_\lambda(z_k)$.
Finally, the evaluation homomorphism $\rho_0$ \eqref{defrho0}
transforms $P_\lambda^{(1/g)}(\mathbf{x})$ into the normalisation
constant $c_\lambda$ \eqref{def-clambda}.
For the evaluation of a Jack polynomial at
$\mathbf{x}=(1,\dots,1)$
see \cite{Stan},  \cite[\S 6, (6.11${}^\prime$)]{Macd}
or \cite{Oko}.
\endproof

\begin{theo}\label{decompS}
The operator
$\mathcal{S}_n:\Cbbd[\mathbf{x}]^{S_n}\mapsto\Cbbd[\mathbf{z}]^{S_n}$,
corresponding to the evaluation homomorphism \eqref{defrho0},
can be decomposed into the product of $n$ operators $\mathcal{A}_k$,
\begin{equation}
\Scal_n=\Acal_1\cdots \Acal_n,
\label{S=AAA}
\end{equation}
\begin{equation}
\Acal_k:\Cbbd[x_1,\dots,x_k]^{S_k}
\mapsto\Cbbd[x_1,\dots,x_{k-1}]^{S_{k-1}}
\otimes\Cbbd[z_k],
\qquad k=1,\dots,n.
\label{actAk}
\end{equation}
In particular, $\Acal_1:\Cbbd[x_1]\mapsto\Cbbd[z_1]$.
The operators $\Acal_k$ in \eqref{S=AAA} are assumed to act trivially
(as unit operators) on the variables $\mathbf{z}$.

Furthermore, the operators $\Acal_k$
act on the {\em restricted polynomials}
$P_\lambda^{(1/g)}(x_1,\allowbreak\dots,\allowbreak x_k,\allowbreak1,\dots,1)$ as follows:
\begin{equation}
\Acal_k:P_\lambda^{(1/g)}(x_1,\dots,x_k,1,\dots,1)
 \mapsto P_\lambda^{(1/g)}(x_1,\dots,x_{k-1},1,\dots,1)\,q_\lambda(z_k).
\label{AP-Pq}
\end{equation}
\end{theo}

The factorisation \eqref{S=AAA}--\eqref{AP-Pq}
possesses the property of dimension reduction $k\mapsto k-1$. 
Such factorisations
have been introduced in \cite{Kuzn} when studying
the inverse problem for integrable models of the 
$sl(2)$ class and have been given the general
name {\it factorised separation chains}.

Note that the formula \eqref{AP-Pq} is consistent with Theorem \ref{decompSx}.
Indeed, applying the chain \eqref{S=AAA} of operators $\Acal_k$ to the
polynomial $P_\lambda^{(1/g)}(x_1,\dots,x_n)$ and using
\eqref{AP-Pq} we reproduce formula \eqref{SP=qqqP}.
The existence of a factorization of the form \eqref{S=AAA} 
is by no means trivial.
One cannot simply define $\Acal_k$ by \eqref{AP-Pq}
because the restricted
Jack polynomials $P_\lambda^{(1/g)}(x_1,\dots,x_k,1,\dots,1)$, as polynomials
of $k$ variables labelled by $n$-dimensional vector $\lambda$,
are not linearly independent.

To prove Theorem \ref{decompS}
we shall describe a recursive procedure for constructing
the operators $\Acal_k$ based on the defining formula
$\Scal_n=\rho_0Q_{z_1}\cdots Q_{z_n}$. The idea is to
pull the homomorphism $\rho_0$ through the sequence
of $Q_{z_k}$'s transforming at each step $Q_{z_k}$ into $\Acal_k$
and producing a new homomorphism $\rho_k$.
Let $\id_{[j,k]}$ for $j\le k$ denote the identity operator in
$\Cbbd[x_j,\dots,x_k]^{S_{k-j+1}}$.

\begin{prop}\label{recurrA}
  There exist a unique sequence of operators $\Acal_k$,
\begin{equation*}
\Acal_k:\Cbbd[x_1,\dots,x_k]^{S_k}
\mapsto\Cbbd[x_1,\dots,x_{k-1}]^{S_{k-1}}
\otimes\Cbbd[z_k],
\qquad k=1,\dots,n,
\end{equation*}
and a unique sequence of homomorphisms $\rho_k$,
\begin{equation*}
\rho_k:\Cbbd[x_{k+1},\dots,x_n]^{S_{n-k}} \mapsto \Cbbd,
\qquad k=0,\dots,n
\end{equation*}
such that:

{\rm(a)} $\rho_0:\Cbbd[\mathbf{x}]^{S_n}\mapsto\Cbbd$
coincides with the evaluation \eqref{defrho0},

{\rm(b)} the operators $\Acal_k$ are normalised by the condition
$\Acal_k:1\mapsto1$,

{\rm(c)} the following relation holds for $k=0,\dots,n-1$:
\begin{equation}
    (\id_{[1,k]}\otimes\rho_{k}) Q_{z_{k+1}}
    =\Acal_{k+1}(\id_{[1,k+1]}\otimes\rho_{k+1}).
\label{rhoQ=Arho}
\end{equation}
The homomorphisms $\rho_k$ defined in this way act on polynomials
$f\in\Cbbd[x_{k+1},\dots,\allowbreak x_n]^{S_{n-k}}$ by restricting
their $n-k$ arguments to $1$:
\begin{equation}
   \rho_k: f(x_{k+1},\dots,x_n) \mapsto f(1,\dots,1), \qquad
  k=0,\dots,n-1.
\label{def-rho_k}
\end{equation}
For $k=n$ we have simply $\rho_n=1$.
\end{prop}

{\bf Proof} of Theorem \ref{decompS} follows then immediately.
Applying recursively formula \eqref{rhoQ=Arho} to
$\Scal_n=\rho_0 Q_{z_1}Q_{z_2}\cdots Q_{z_n}$ we get
$\Scal_n=\Acal_1(\id_{[1,1]}\otimes\rho_1)Q_{z_2}\cdots Q_{z_n}$ and so on,
until we arrive at \eqref{S=AAA}. Applying the identity
\eqref{rhoQ=Arho} to a Jack polynomial $P_\lambda^{(1/g)}(\mathbf{x})$
and using \eqref{def-rho_k} we get \eqref{AP-Pq}.
\endproof

{\bf Proof} of Proposition \ref{recurrA} is given by induction.
For $\rho_0$ the formula \eqref{def-rho_k} holds by assumption.
Assuming that for some $k\in[0,n-1]$
the operators $\Acal_1,\dots,\Acal_k$ and homomorphisms
$\rho_0,\dots,\rho_k$ have been already constructed, let us consider
the product $(\id_{[1,k]}\otimes\rho_k)Q_{z_{k+1}}$
and show that it factorises uniquely as
$\Acal_{k+1}(\id_{[1,k+1]}\otimes\allowbreak\rho_{k+1})$.
Note that the operation $(\id_{[1,k]}\otimes\rho_k)$ produces
{} from a polynomial $f(\mathbf{x})$ of $n$ variables
the polynomial $f(x_1,\dots,x_k,1,\dots,1)$ of $k$ variables.
As in Section \ref{q}, we will distinguish the variables $\mathbf{y}$
of the target space for the operator $Q_{z_{k+1}}$ from the variables
$\mathbf{x}$ of the source space. Given a polynomial
$f\in\Cbbd[\mathbf{x}]^{S_n}$, we have thus to take
$Q_{z_{k+1}}f\in\Cbbd[\mathbf{y}]^{S_n}$ and to set $y_{k+1}=\ldots=y_n=1$.
We can analyse the resulting expression in two ways based on two
descriptions of $Q_{z_{k+1}}$ given in Section \ref{q}: the one in terms of an
integral operator and a pure algebraic one.

Let us apply first the combinatorial prescription for evaluating
$Q_{z_{k+1}}$ on the elementary symmetric polynomials
$e_j(\mathbf{x})$ given in Remark \ref{algorithm},
formulae \eqref{qop36} and \eqref{subs-eta}.
Recalling formula \eqref{defw},
we introduce the generating functions
$w_x(t)=\prod_{j=1}^n(t-x_j)$ for $e_j(\mathbf{x})$ and, respectively,
$w_y(t)=\prod_{j=1}^n(t-y_j)$ for $e_j(\mathbf{y})$.
Consider now the relation \eqref{qop13},
\begin{equation}
 w_x(t)=w_y(t)-\sum_{i=1}^n\eta_i y_i\,\frac{w_y(t)}{t-y_i}\,,
\label{wxwy}
\end{equation}
and set $y_{k+1}=\ldots=y_n=1$. The right hand side  of \eqref{wxwy} is then divisible
by $(t-1)^{n-k-1}$, so is the left hand side, therefore $n-k-1$ of the variables
$x_j$ are forced to take value $1$. Due to the symmetry, the remaining
free variables can be chosen as $x_1,\dots,x_{k+1}$. Note that
{} from the above the statement of the Proposition follows immediately:
we have effectively transformed the operator expression
$(\id_{[1,k]}\otimes\rho_k)Q_{z_{k+1}}$ to the desired form
$\Acal_{k+1}(\id_{[1,k+1]}\otimes\rho_{k+1})$ with some
operator $\Acal_{k+1}$ acting on the variables $x_1,\dots,x_{k+1}$.

We can go even further and find how $\Acal_{k+1}$ acts on polynomials
in $x_1,\dots,x_{k+1}$.
Let 
\begin{equation*}
w_x(t)=(t-1)^{n-k-1}\widetilde{w}_x(t),
\qquad \widetilde{w}_x(t)\equiv\prod_{j=1}^{k+1}(t-x_j)
\end{equation*}
and, respectively,
\begin{equation*}
w_y(t)=(t-1)^{n-k-1}\widetilde{w}_y(t),
\qquad \widetilde{w}_y(t)\equiv\prod_{j=0}^{k}(t-y_j),
\qquad y_0\equiv1.
\end{equation*}
The relation \eqref{wxwy} transforms then into
\begin{equation}
 \widetilde{w}_x(t)=\widetilde{w}_y(t)-\sum_{i=0}^k\eta_i y_i\,
 \frac{\widetilde{w}_y(t)}{t-y_i}\,,
\qquad \eta_0\equiv \eta_{k+1}+\ldots+\eta_n.
\label{wxwy-tilde}
\end{equation}
Expanding both sides of \eqref{wxwy-tilde} in $t$ we get
the following prescription for calculating $\Acal_{k+1}f$:
\begin{enumerate}
\item express $f\in\Cbbd[x_1,\dots,x_{k+1}]^{S_{k+1}}$ 
as a polynomial in $\{e_j(x_1,\dots,x_{k+1})\}$, $j=1,\dots,k+1$;
\item substitute for $e_j(x_1,\dots,x_{k+1})$ the expressions 
\begin{equation}
  e_j(x_1,\dots,x_{k+1})=(1+\eta_0+\ldots+\eta_{j-1})y_0\cdots y_{j-1}
+\mbox{\rm permutations}
\label{ejxejy}
\end{equation}
(the distinct permutations are taken simultaneously in
$y_j$ and $\eta_j$);
\item expand the resulting polynomial in monomials
$\eta_0^{m_0}\eta_1^{m_1}\cdots\eta_{k}^{m_k}$, 
$\mathbf{m}=(m_0,\allowbreak m_1,\dots,m_k)\in
\Nbbd^{k+1}$;
\item replace each monomial in $\eta_j$ according to the rule
\begin{equation}
  \eta_0^{m_0}\eta_1^{m_1}\cdots\eta_{k}^{m_k}
  \mapsto\frac{(z_{k+1}-1)^{m_0+\ldots+m_k}((n-k)g)_{m_0}}{(ng)_{m_0+\ldots+m_k}}
  \prod_{j=1}^k (g)_{m_i}.
\label{etaynew}
\end{equation}
\end{enumerate}
The last formula is obtained 
by expanding 
$(\eta_{k+1}+\ldots+\eta_n)^{m_0}$
in the left hand side of \eqref{etaynew},
evaluating the result by formula \eqref{subs-eta}
and using the following
combinatorial identity for $s\equiv n-k$, $m\equiv m_0$:
\begin{equation}
  \sum_{j_1+\ldots+j_s=m} 
  \frac{m!}{j_1!\cdots j_s!} \prod_{\alpha=1}^s (g)_{j_\alpha}
= (sg)_m,
\label{comb-id}
\end{equation}
which is easily proved by expanding the identity
$[(1-t)^{-g}]^s=(1-t)^{-sg}$ in powers of $t$.

An alternative proof of Proposition \ref{recurrA} is based on the
formula \eqref{qop7}, presenting $Q_{z_{k+1}}$ as an integral operator.
Let us take the expression \eqref{qop7} for $[Q_{z_{k+1}}f](\mathbf{y})$
and restrict $n-k$ of the arguments $y_j$ to the value $1$.
Since $[Q_{z_{k+1}}f](\mathbf{y})$ is a symmetric polynomial, it does not matter
which of $y_j$ we choose to fix. Because of the inequalities in the 
definition of the domain $\Omega_\mathbf{x}$ \eqref{liov6},
\begin{equation*}
0<y_1<x_1<y_2<\ldots<y_n<x_n,
\end{equation*}
the natural choice is to take the limit
$y_1,y_2,\dots,y_{n-k}\rightarrow\widetilde{y}_0\equiv 1$.
The variables $x_1,\dots,x_{n-k-1}$ are sandwiched between $y$'s
and, therefore, forced to tend to $1$ as well.
As $y_1,\dots,y_{n-k}$ tend to $y_0\equiv1$, the right hand side
 of \eqref{qop7}
exhibits an uncertainty $0/0$. To resolve the uncertainty, let us set
$y_j=1+\varepsilon v_j$, $j=1,\dots,n-k$, and
$x_j=1+\varepsilon u_j$, $j=1,\dots,n-k-1$,
for some $\varepsilon>0$ and $u_j$, $v_j$ satisfying
\begin{equation*}
  0<v_1<u_1<\ldots<u_{n-k-1}<v_{n-k},
\end{equation*}
and take the limit $\varepsilon\rightarrow0$.
Let us also renumber the
remaining arguments as $y_{n-k+j}=\widetilde{y}_j$,
$j=1,\dots,k$, and $x_{n-k+j-1}=\widetilde{x}_j$,
$j=1,\dots,k+1$.

The components of the integrand of \eqref{qop7} are then transformed as
follows:
\begin{eqnarray*}
 \xi_i &\rightarrow& \hat{\xi}_i\widetilde{\xi}_0,
\qquad i=1,\dots,n-k, \\
 \xi_{n-k+i} &\rightarrow& \widetilde{\xi}_i,
\qquad i=1,\ldots k,
\end{eqnarray*}
where
\begin{eqnarray*}
 \hat{\xi}_i&=&
 \frac{\prod_{j=1}^{n-k-1}(u_j-v_i)}%
{\prod^{n-k}_{\genfrac{}{}{0pt}{}{j=1}{j\neq i}} (v_j-v_i)}\,,
\qquad i=1,\dots,n-k, \\
 \widetilde{\xi}_i &=&
 \frac{\prod_{j=1}^{k+1}(\widetilde{x}_j-\widetilde{y}_i)}%
{\prod^{k}_{\genfrac{}{}{0pt}{}{j=0}{j\neq i}}
(\widetilde{y}_j-\widetilde{y}_i)}\,,
\qquad i=0,\dots,k;
\end{eqnarray*}
\begin{equation*}
 \frac{\Delta_n(\mathbf{x})}{\Delta_n(\mathbf{y})}\sim
 \varepsilon^{-n+k+1}\widetilde{\xi}_0^{n-k-1}
 \frac{\Delta_{n-k-1}(\mathbf{u})}{\Delta_{n-k}(\mathbf{v})}\,
 \frac{\Delta_{k+1}(\boldsymbol{\widetilde{x}})}%
{\Delta_{k+1}(\boldsymbol{\widetilde{y}})}
\end{equation*}
and $\boldsymbol{\widetilde{x}}=(\widetilde{x}_1,\ldots,
\widetilde{x}_{k+1})$, 
$\boldsymbol{\widetilde{y}}=(\widetilde{y}_0,\ldots,\widetilde{y}_k)$.

Since $dx_1\cdots dx_{n-k-1}\sim\varepsilon^{n-k-1}du_1\cdots du_{n-k-1}$,
the factors $\varepsilon$ cancel completely. Integration in the variables
$u_j$ produces a constant factor which is evaluated 
with the help of Proposition \ref{prop7}:
\begin{equation}
  \frac{1}{\Delta_{n-k}(\mathbf{v})}
  \int_{v_1}^{v_2}du_1\cdots\int_{v_{n-k-1}}^{v_{n-k}}du_{n-k-1}\,
\Delta_{n-k-1}(\mathbf{u})\,\left[
\prod_{j=1}^{n-k}\hat{\xi}_j^{g-1}\right]
 =\frac{\Gamma(g)^{n-k}}{\Gamma((n-k)g)}\,.
\end{equation}
The expression $(\id_{[1,k]}\otimes\rho_k)Q_{z_{k+1}}$ is presented in the
desired form $\Acal_{k+1}(\id_{[1,k+1]}\allowbreak\otimes\rho_{k+1})$ with the
operator $\Acal_{k+1}$ defined by the integral
\begin{multline}
 [\Acal_{k+1}f](\boldsymbol{\widetilde{y}})=
  \frac{\Gamma(ng)}{\Gamma((n-k)g)\Gamma(g)^k}\,
  \frac{(z_{k+1}-1)^{1-ng}}{\Delta_{k+1}(\boldsymbol{\widetilde{y}})}\,
\left[\prod_{j=1}^k \widetilde{y}_j^{1-g}  \right]\\
\times \int_{\widetilde{\Omega}_\mathbf{x}} d\boldsymbol{\widetilde{x}}\,
\Delta_{k+1}(\boldsymbol{\widetilde{x}})\,
\delta\left(\widetilde{x}_1\cdots\widetilde{x}_{k+1}-
z_{k+1}\widetilde{y}_1\cdots\widetilde{y}_{k}\right)\,
\widetilde{\xi}_0^{(n-k)g-1}
\left[\prod_{j=1}^k \widetilde{\xi}_j^{g-1}\right]\, f(\boldsymbol{\widetilde{x}}),
\label{Aint}
\end{multline}
where the integration domain $\widetilde{\Omega}_\mathbf{x}$ is described by the
inequalities
\begin{equation*}
0<1\equiv\widetilde{y}_0<\widetilde{x}_1<\widetilde{y}_1
<\ldots<\widetilde{y}_k<\widetilde{x}_{k+1}<\infty.
\end{equation*}

Using the Dirichlet-Liouville integral \eqref{liov2}, one can analyse the
integral operator \eqref{Aint} in the same manner as the operator $Q_z$
has been analysed in section 4. 
In particular, one can prove an analogue of Proposition
\ref{prop2}: $\Acal_{k+1}$ sends symmetric polynomials into symmetric
polynomials (with the dimension reduced by one)
 and its action on the elementary symmetric polynomials
$e_j$ is described by formulae \eqref{ejxejy} and \eqref{etaynew}.
\endproof

Notice that, for $k=0$, the formula \eqref{Aint} produces 
$\Acal_1:f(x_1)\mapsto f(z_1)$, whereas setting $k=1$ in \eqref{AP-Pq}
we get
\begin{equation*}
\Acal_1:P_\lambda^{(1/g)}(x_1,1,\dots,1)
 \mapsto P_\lambda^{(1/g)}(1,\dots,1)\,q_\lambda(z_1)
 =c_\lambda q_\lambda(z_1).
\end{equation*}
Combining these two observations together, we arrive at the
following remarkable expression for the separated polynomial
$q_\lambda(z)$ in terms of a restricted Jack polynomial:
\begin{equation}
 q_\lambda(z)=c_\lambda^{-1} \, P_\lambda^{(1/g)}(z,1,\dots,1).
\label{amazing}
\end{equation}

\section{Integral representation}\label{intrepr}
\noindent
In this section we shall apply our $Q$-operator to construct an integral
representation for Jack polynomials. In principle, such representation
could be obtained by inverting the separating operator $\Scal_n$,
described in Theorem \ref{decompS}.
However, the construction of $\Scal_n^{-1}$ is a difficult 
and yet unsolved problem.
Fortunately, there is another, more direct approach using the operator
$Q_z$ at $z=0$.

Applying Proposition \ref{commQe} to the case $z=0$ we conclude
that $Q_0$ nullifies the
ideal $e_n(\mathbf{x})\Cbbd[\mathbf{x}]^{S_n}$ and therefore can be
canonically defined on the factor space
$\Fcal=\Cbbd[\mathbf{x}]^{S_n}/e_n(\mathbf{x})\Cbbd[\mathbf{x}]^{S_n}$.
For the rest of this section we denote $\mathbf{x}=(x_1,\dots,x_n)$,
$\mathbf{y}=(y_1,\dots,y_n)$ and $\mathbf{x}^\prime=(x_1,\dots,x_{n-1})$.
Let $\Pcal$ be the projection operator
$\Pcal:\Cbbd[\mathbf{x}]^{S_n}\mapsto\Cbbd[\mathbf{x}^\prime]^{S_{n-1}}:
p\,(x_1,\dots,x_n)\mapsto p\,(x_1,\dots,x_{n-1},0)$.
Note that $\Pcal:e_i(\mathbf{x})\mapsto e_i(\mathbf{x}^\prime)$ for
$i=1,\dots,n-1$ and $\Pcal:e_n(\mathbf{x})\mapsto0$.
Since the products $e_1^{k_1}\cdots e_{n-1}^{k_{n-1}}$ form a basis in
$\Fcal$, the pro\-jec\-tor $\Pcal$ provides a natural isomorphism
$\Fcal\simeq\Cbbd[\mathbf{x}^\prime]^{S_{n-1}}$. We have thus come to
the following conclusion.

\begin{prop}\label{qprime}
There exists a unique operator
$Q_0^\prime:\Cbbd[\mathbf{x}^\prime]^{S_{n-1}}\mapsto\Cbbd[\mathbf{y}]^{S_n}$
such that
\begin{equation}
  Q_0=Q_0^\prime\Pcal.
\label{def-qprime}
\end{equation}
\end{prop}
Formula \eqref{def-qprime} provides a direct way to constructing
a Jack polynomial in $n$ variables from a Jack polynomial in $n-1$
variables. For any partition $\lambda=(\lambda_1,\dots,\lambda_n)$
of length $n$ define two partitions:
$\lambda^\flat=(\lambda_{1,n},\lambda_{2,n},\dots,\lambda_{n-1,n},0)$
of length $n$
and $\lambda^\natural=(\lambda_{1,n},\lambda_{2,n},\dots,\lambda_{n-1,n})$
of length $n-1$.
Recall two important properties of Jack polynomials
(see \cite{Macd}): homogeneity,
\begin{equation}
P^{(1/g)}_{\l_1,\dots,\l_n}(x_1,\dots,x_n)=(x_1\cdots x_n)^{\l_n}
P^{(1/g)}_{\l_{1,n},\dots,\l_{n-1,n},0}(x_1,\dots,x_n),
\label{integr9}
\end{equation}
and restriction to $x_n=0$,
\begin{equation}
P^{(1/g)}_{\l_{1,n},\dots,\l_{n-1,n},0}(x_1,\dots,x_{n-1},0)=
P^{(1/g)}_{\l_{1,n},\dots,\l_{n-1,n}}(x_1,\dots,x_{n-1}).
\label{integr10}
\end{equation}
Note that the equations \eqref{integr9} and \eqref{integr10} can be abbreviated
as $P^{(1/g)}_{\lambda}(\mathbf{x})
=\left[e_n(\mathbf{x})\right]^{\l_n}\allowbreak P^{(1/g)}_{\lambda^\flat}(\mathbf{x})$ and, respectively,
$\Pcal:P^{(1/g)}_{\lambda^\flat}(\mathbf{x})
\mapsto P^{(1/g)}_{\lambda^\natural}(\mathbf{x}^\prime)$.

By Theorem \ref{qop41}, the eigenvalue of $Q_0$ on
$P^{(1/g)}_{\lambda}(\mathbf{x})$ is given by $q_\lambda(0)$.
It follows from Theorem~\ref{theorem1} that, for small $z$, 
$f_\l(z)=z^{\l_n}+O(z^{\l_n+1})$ and then from
Theorem~\ref{theorem2} that
\begin{equation}
q_\l(z)=\beta_\l^{-1}z^{\l_n}(1+O(z)), \qquad z\rightarrow0,
\label{integr11}
\end{equation}
where $\beta_\l$ is defined by (\ref{conn19}). In particular,
$q_{\lambda^\flat}(0)=\beta_\lambda^{-1}$.

Applying the operator equality \eqref{def-qprime} to the polynomial
$P^{(1/g)}_{\lambda^\flat}(\mathbf{x})$ and using \eqref{integr10}
and $q_{\lambda^\flat}(0)=\beta_\lambda^{-1}$
we obtain
$\beta_\lambda^{-1}P^{(1/g)}_{\lambda^\flat}(\mathbf{x})
=[Q_0^\prime P^{(1/g)}_{\lambda^\natural}](\mathbf{x})$.
Finally, we use \eqref{integr9} and arrive at the formula
\begin{equation}
   P^{(1/g)}_{\lambda}(\mathbf{x})=
   \beta_\lambda \left[e_n(\mathbf{x})\right]^{\lambda_n}
   [Q_0^\prime P^{(1/g)}_{\lambda^\natural}](\mathbf{x}),
\label{reconstrP}
\end{equation}
expressing a Jack polynomial in $n$ variables in terms of
a Jack polynomial in $n-1$ variables.

To give the formula \eqref{reconstrP} more flesh, we need
to realise $Q_0$ as an integral operator. However,
formula \eqref{qop7}, which we used to define $Q_z$, is valid only for
$z>1$, so that the case $z=0$ requires a special consideration.

\begin{prop}\label{prop7}
Given a set of real parameters
$\mathbf{y}=(y_1,\dots,y_n)$ satisfying inequalities \eqref{liov2a},
define the domain $\Omega_\mathbf{x}^\prime\subset \Rbbd^{n-1}$ by the inequalities
\begin{equation}
\Omega_\mathbf{x}^\prime=
\{\mathbf{x}^\prime\in \Rbbd^{n-1} \mid
 0< y_1< x_1< y_2<\ldots< x_{n-1}< y_n<\infty\}.
\label{Omega-prime}
\end{equation}
Then for $\Re\a_i>0$
\begin{equation}
\int_{\Omega_\mathbf{x}^\prime} d\mathbf{x}^\prime\,
\frac{\Delta_{n-1}(\mathbf{x}^\prime)}{\Delta_n(\mathbf{y})}\,
\left[\prod_{i=1}^{n} \hat\xi_i^{\a_i-1}\right]=
\frac{\Gamma(\a_1)\cdots\Gamma(\a_n)}{\Gamma(\abs{\boldsymbol{\a}})}\,,
\label{integr3}
\end{equation}
where $\hat\xi_i>0$ are defined as
\begin{equation}
\hat{\xi}_i=\frac{\prod_{k=1}^{n-1}(x_k-y_i)}%
{\prod^n_{\genfrac{}{}{0pt}{}{k=1}{k\neq i}} (y_k-y_i)}\,,
\quad i=1,\dots,n,\quad
\sum_{i=1}^n\hat\xi_i=1.
\label{integr2}
\end{equation}
\end{prop}

The {\bf proof} follows closely that of Theorem \ref{dirichlet1}.
Introducing functions
$\hat{w}_x(t)=\prod_{i=1}^{n-1}(t-x_i)$ and
$w_y(t)=\prod_{i=1}^{n}(t-y_i)$ we construct the generating function
$\hat{R}(t)=\hat{w}_x(t)/w_y(t)=\sum_{i=1}^n \hat\xi_i/(t-y_i)$
for the change of variables
$(x_1,\dots,x_{n-1})\mapsto(\hat\xi_1,\dots,\hat\xi_{n-1})$.
The argument, based on studying $\hat{R}(t)$, shows that the domain
$\Omega_\mathbf{x}^\prime$ is mapped bijectively onto $V_{\boldsymbol{\xi}}^\prime$,
\begin{equation}
V_{\boldsymbol{\xi}}^\prime=\{\hat{\boldsymbol{\xi}^\prime}\in\Rbbd^{n-1} \mid
\hat\xi_i>0,\; i=1,\dots,n-1, \quad
\hat{\xi_1}+\ldots+\hat\xi_{n-1}<1 \}.
\label{def-Vxiprime}
\end{equation}
The corresponding Jacobian is given by
$ \det \left[\partial\xi_j/\partial x_k \right]_{i,j=1}^n=
 \Delta_{n-1}(\mathbf{x}^\prime)/\Delta_n(\mathbf{y})$.

Making the change of variables
$\mathbf{x}^\prime\mapsto\boldsymbol{\hat\xi^\prime}$
in the integral (\ref{integr3}) we come
exactly to the Dirichlet-Liouville integral (\ref{liov2})
with the only difference that the integration
variable $\hat\xi_n$ is eliminated by resolving explicitly the
constraint $\abs{\hat{\boldsymbol{\xi}}}=1$.
\endproof

\begin{theo}\label{Q0}
The integral operator $\Qbbd$ defined by
\begin{equation}
[\Qbbd\, p\,](\mathbf{y})=
\frac{\Gamma(ng)}{\Gamma(g)^n}
\int_{\Omega_\mathbf{x}^\prime} d\mathbf{x}^\prime\,
\frac{\Delta_{n-1}(\mathbf{x}^\prime)}{\Delta_n(\mathbf{y})}\,
\left[\prod_{i=1}^{n} \hat\xi_i^{g-1}\right]\,
p\,(\mathbf{x}^\prime)
\label{integr7}
\end{equation}
sends $\Cbbd[\mathbf{x}^\prime]^{S_{n-1}}$ into
$\Cbbd[\mathbf{y}]^{S_{n}}$ and coincides with the operator
$Q_0^\prime$ defined by \eqref{def-qprime}.
\end{theo}

{\bf Proof.} Repeating almost word by word the proof of Proposition \ref{prop2}
and using formula \eqref{integr3}, we find the analogue of the
formula \eqref{qop11} $\forall \,\mathbf{k}=(k_1,\ldots,k_n)\in\Nbbd^n$,
\begin{equation}
  \Qbbd: \hat{\boldsymbol{\xi}}^{\mathbf{k}}
   \mapsto
    \frac{(g)_{k_1}\cdots(g)_{k_n}}{(ng)_{\abs{\mathbf{k}}}}\,.
\label{qprimeint}
\end{equation}
The action of $\Qbbd$ on the polynomials $e_i(\mathbf{x}^\prime)$
is obtained then from the analogue of the formula \eqref{qop14},
\begin{equation}
e_i(\mathbf{x}^\prime)=
   \sum_{j=1}^n \hat\xi_j \left[e_{i-1}(\mathbf{y})\right]_{y_j=0},
\label{qop14a}
\end{equation}
which is obtained, in turn, by expanding the equality
$\hat{w}_x(t)=\sum_{i=1}^n\hat\xi_i w_y(t)/(t-y_i)$ in powers of $t$.

It now remains to identify the action of $\Qbbd$ and $Q_0^\prime$
on $\Cbbd[\mathbf{x}^\prime]^{S_{n-1}}$.
Setting $x_n=0$ in the expressions \eqref{liov5} for $\xi_i$ and comparing the
result with the expressions \eqref{integr2}, we obtain
\begin{equation}
  \Pcal\xi_i=-y_i\hat\xi_i.
\label{xixi}
\end{equation}
Setting $z=0$ in the right hand side
of \eqref{qop11} and using \eqref{xixi} we obtain
exactly the right hand side of \eqref{qprimeint}, which proves that
$Q_0=\Qbbd\,\Pcal$. 
\endproof

We can regard the formula \eqref{reconstrP} as an integral relation,
which reduces a Jack polynomial in $n$ variables to a Jack polynomial
in $n-1$ variables\footnote{or, vice versa, builds up the number 
of variables, starting from the unit function}. 
Starting with $P_\l^{(1/g)}(x_1,\ldots,x_n)$ and
having iterated the formula \eqref{reconstrP} $n-1$ times, 
we go down to the $P_{\l_1}^{(1/g)}(x_1)=x_1^{\l_1}$ and,
as a result, get an explicit integral representation 
for a Jack polynomial in $n$ variables.
To describe it explicitly,
we introduce a triangular matrix $x_{i,j}$, $1\le i\le j\le n$
of variables.
At the first step, the variables $x_{i,n}$, $i=1,\dots,n$, are the external ones and
 $x_{i,n-1}$, $i=1,\dots,n-1$, are the integration
variables. Then the equation \eqref{reconstrP} takes the following form:
\begin{multline}
 P_{\l}^{(1/g)}(x_{1,n},\dots,x_{n,n})=
\left[\prod_{i=1}^n x_{i,n}^{\l_n}\right]\,
\left[\prod_{i=1}^{n-1}\frac{((n-i+1)g)_{\l_{i,n}}}{((n-i)g)_{\l_{i,n}}}\right] \\
 \times\int_{x_{1,n}}^{x_{2,n}}dx_{1,n-1}\cdots\int_{x_{n-1,n}}^{x_{n,n}}
dx_{n-1,n-1} \\
\times
\frac{\prod_{i<j}(x_{j,n-1}-x_{i,n-1})}{\prod_{i<j}(x_{j,n}-x_{i,n})} \,
\left[\frac{\ds\prod_{i=1}^{n}
\prod_{k=1}^{n-1}(x_{i,n}-x_{k,n-1})}
{\ds\prod_{k\neq i}(x_{i,n}-x_{k,n})}
\right]^{g-1} \,
\\
\times P_{\l_{1,n},\dots,\l_{n-1,n}}^{(1/g)}(x_{1,n-1},\dots,
x_{n-1,n-1}).\qquad \qquad \qquad
\label{integr12}
\end{multline}
Iterating this formula, we obtain the integral representation for
Jack polynomials
\begin{multline}
P_{\l}^{(1/g)}(x_{1,n},\dots,x_{n,n})=
\frac{\ds\prod_{i=1}^nx_{i,n}^{\l_n}}
{\ds\prod_{i<j}(x_{j,n}-x_{i,n})}\,
\left[\prod_{j=2}^n\prod_{i=1}^{j-1}\frac{((j+1-i)g)_{\l_{i,j}}}
{((j-i)g)_{\l_{i,j}}} \right]\\
 \times\Biggl(\prod_{{\phantom{aa}
\stackrel{\longleftarrow}{ j={n-1,1}}}}
\biggl(\prod_{\phantom{aa}\stackrel{\longleftarrow}{{i=1,j}}}
\int_{x_{i,j+1}}^{x_{i+1,j+1}}dx_{i,j}\biggr)\Biggr)
\Biggl\{\prod_{j=2}^n\prod_{i=1}^j\biggl[\>
\frac{\prod_{k=1}^{j-1}(x_{i,j}-x_{k,j-1})}
{\prod_{1=k\neq i}^j
{\ds(x_{i,j}-x_{k,j})}}\biggr]\Biggr\}^{g-1} \\
 \times
\prod_{j=1}^{n-1}\prod_{i=1}^j x_{i,j}^{\l_{j,j+1}}.
\label{repr1}
\end{multline}
It is implied in (\ref{repr1}) that the integrals $\int dx_{i,j}$
are ordered corresponding to the order of indices in the products.

As mentioned before, formula \eqref{integr12} first appeared
in \cite{Oko}, where its proof was
based on the fact that Jack polynomials can be obtained
as a limit from the shifted Jack polynomials.
One can regard the integral equation \eqref{qop42} in section 4
as the one-parameter generalisation of \eqref{integr12}.

\section{Concluding remarks}\label{discuss}
\noindent
Our construction of the $Q$-operator relies on
the properties of Jack polynomials. Nevertheless, several of the 
obtained results are of importance for the general theory
of quantum integrable models. It concerns, first of all,
the {\it factorised separation chain}, constructed in Section \ref{restricted}.
Similar factorisations of separating operators were already observed for
other integrable models \cite{Kuzn} and, apparently, they are manifestations
of some general pattern, which is yet to be fully understood.

Notice also the intriguing role of the {\it restricted Jack polynomials} revealed
in formulae like \eqref{AP-Pq} or \eqref{amazing}.
It would be interesting to find  analogous formulae 
for other integrable models.
The properties of restricted Jack polynomials themselves deserve a
further study.

Finally, we should mention another integral representation known
for Jack polynomials \cite{Awata}, where, in contrast with \eqref{repr1},
the number of integrations
grows with the weight of $\l$.

\section{Acknowledgements}\label{ack}
\noindent
VBK and EKS would like to acknowledge the support from
the Isaac Newton Institute (Cambridge, UK) during the programme
``Symmetric Functions and Macdonald Polynomials'' (2001).

VVM acknowledges the support from Australian Research Council,
including support for research visits between ANU and  University of Leeds. 
VBK is supported  by EPSRC.


\end{document}